\documentclass[12pt]{amsart}
\usepackage{amsmath,amsthm,amssymb,enumerate,mathscinet,mathtools}
\usepackage{fullpage}
\usepackage[svgnames,dvipsnames]{xcolor}
\usepackage{graphicx,xcolor,pgfplots}
\usepackage{verbatim}
\usepackage{mathrsfs}
\usepackage{caption}
\usepackage{hyperref}
\usepackage{enumitem}
\usepackage{blindtext}
\usepackage{multicol}
\usetikzlibrary{fit,calc,positioning,decorations.pathreplacing,matrix}
\newtheorem{theorem}{Theorem}[section]
\newtheorem{fact}[theorem]{Fact}
\newtheorem{lemma}[theorem]{Lemma}

\newtheorem{proposition}[theorem]{Proposition}
\newtheorem{corollary}[theorem]{Corollary}

\newtheorem{definition}[theorem]{Definition}

\newtheorem{thm}[theorem]{Theorem}

\theoremstyle{definition}

\newtheorem{remark}[theorem]{Remark}

\numberwithin{equation}{section}

\def\C{\mathcal{C}}

\def\N{\mathbb{N}}

\def\D{{\color{blue} D}}
\def\U{{\color{blue} U}}
\def\Fb{{\color{blue} F_1}}
\def\Fg{{\color{Green} F_2}}
\def\F{\mathcal{F}}
\def\FF{\mathbb{F}}
\def\u{(1,1)}
\def\d{(1,-1)}

\makeatletter
\newcommand*{\rom}[1]{\expandafter\@slowromancap\romannumeral #1@}

\makeatother
\newcommand{\seqnum}[1]{\href{http://oeis.org/#1}{\underline{#1}}}

\def\COMMENT#1{}
\let\COMMENT=\footnote

\allowdisplaybreaks

\title{Intervals in a family of Fibonacci lattices}

\date{\today}

\author{}

\thanks{}

\begin{document}
\author[J.-L Baril]{Jean-Luc Baril}

\author[N. Hassler]{Nathana\"{e}l Hassler}
\address{LIB, Universit\'e de Bourgogne Franche-Comt\'e,   B.P. 47 870, 21078, Dijon Cedex, France}
\email{barjl@u-bourgogne.fr, nathanael.hassler@ens-rennes.fr}

\begin{abstract}
We focus on a family of subsets $(\F^p_n)_{p\geq 2}$ of Dyck paths of semilength $n$ that avoid the patterns $DUU$ and $D^{p+1}$, which are enumerated by  the generalized Fibonacci numbers. We endow them with the partial order relation induced by the well-known Stanley lattice, and we prove that all these posets are sublattices of the Stanley lattice.  We provide generating functions for the numbers of linear  and boolean intervals  and we deduce  the Möbius function  for every $p\geq 2$. We  count meet-irreducible elements in $\FF_n^p$ which establishes a surprising link with the edges of the $(n,p)$-Tur\'an graph. We also prove that intervals are in one-to-one correspondence with bicolored Motzkin paths avoiding some patterns, which allows to enumerate intervals for $p=2$. Using a discrete continuity argument ($p\rightarrow \infty$), we present a similar enumerative study in a poset of some Dyck paths of semilength $n$ counted by $2^{n-1}$. Finally, we give bijections that transport the lattice structure on other combinatorial objects, proving that those lattices can be seen as the well-known dominance order on some compositions.
\end{abstract}

\maketitle

\section{Introduction and notation}\label{sec:intro}

The Fibonacci sequence, defined by the recurrence relation $F_n=F_{n-1}+F_{n-2}$ for $n\geq 2$, with initial conditions $F_0=F_1=1$, is one of the most fascinating integer sequences in mathematics. It is probably due to its elegant properties and its appearance in various natural phenomena in different domains as biology, computer science, finance, architecture (see for instance \cite{Knu,Kos}). This classical sequence lays the groundwork for a broader class of sequences: the $p$-generalized Fibonacci sequences defined for every $p\geq 2$ by $$F_n^p=F_{n-1}^p+ F_{n-2}^p+\cdots +F_{n-p}^p$$ with initial conditions $F_i^p=0$ for $i<0$, and $F_0^p=1$ (see \cite{Mil}). These sequences are frequently found in the literature. In number theory, numerous studies introduce new formulas and properties for these numbers, while combinatorial research uncovers new classes of combinatorial objects counted by these numbers and, at times, develops efficient algorithms for their complete generation. For instance, binary words of length $n$ avoiding the pattern  $1^p$ (for a given $p\geq 2$) are counted by the $p$-generalized Fibonacci number $F_{n+1}^p$, and  their exhaustive generation can be obtained in Gray code order with a constant amortized time  algorithm \cite{Vaj}. However, there is a bit less work studying a partial order  on a combinatorial class counted by $F_n$. For instance, Stanley  proved  in \cite{Stanley} that the Young-Fibonacci poset $\texttt{Z}(r)$ (also called Fibonacci $r$-differential poset)  and the $r$-Fibonacci poset $\texttt{Fib}(r)$ are two lattices \cite{KO,Kre,Stan,Stanley,Stfur}. To our knowledge, there is no work that study partial order on a family of combinatorial classes $\mathcal{C}(p)$, $p\geq 0$, where the elements of size $n$ in $\mathcal{C}(p)$ are counted by the generalized Fibonacci numbers $F_n^p$. This is one of the key aims of this study.

 A {\it Dyck path} of semilength $n$, $n\geq 0$, is a lattice path  in the first quarter plane starting at the origin  $(0, 0)$, 
ending at $(2n, 0)$, and never going below the $x$-axis, consisting of up steps
$U = (1, 1)$ and down steps $D=(1,-1)$. We write $\epsilon$ the empty path, that is the only path with semilength $0$. Let $\mathcal{D}$ be the set of all Dyck paths, and $\mathcal{D}_n$ be the set of  those of semilength $n$. For instance, we have $$\mathcal{D}_3=\{UUUDDD, UUDDUD, UUDUDD, UDUDUD, UDUUDD\}.$$ A {\it pattern} in a Dyck path $P$ consists of consecutive steps of $P$. For instance a {\it peak} $UD$  is a pattern that always appears in a nonempty Dyck path, while a {\it valley} $DU$ does not occur in $U^nD^n$, $n\geq 0$. More generally, we will say that a Dyck path $P$ avoids a pattern $\alpha$ when $P$ does not contain the factor $\alpha$. Let $P$ be a Dyck path, we define $\texttt{type}(P)$ as the length of the last descent run  of $P$. For instance, we have $\texttt{type}(UDUUDD)=2$ since the last descent run is $DD$.  

There exist several partial ordering relations on Dyck paths which endows them with lattice structures \cite{Barc,Barp1,Barp2,BKN,BBKN,Chapo,Simi,enum comb, Tam}. Of much interest are probably
the so-called Tamari lattice~\cite{Huan,Tam} obtained with the  rotations on Dyck paths
\cite{Barmo,Berg,Chapo}, and the Stanley lattice \cite{enum comb} obtained with the covering $DU\rightarrow UD$ and where  $P\leq Q$ if and only if $P$ is always below $Q$ when we draw them in the quarter plane.

Throughout this paper, we focus on the family of sets $\F^p$, $p\geq 2$, where the elements are Dyck paths in $\mathcal{D}$ avoiding the patterns $DUU$ and $D^{p+1}$. The set of elements in $\F^p$ having semilength $n$ will be denoted $\F_n^p$. We will also consider the set $\F^\infty$ (resp. $\F_n^\infty$) of Dyck paths (resp. of semilength $n$) avoiding the pattern $DUU$. Note the following inclusions for all $n\geq0$, which will allow us to use discrete continuity arguments:
$$\F_n^2\subseteq\F_n^3\subseteq\cdots\subseteq\F_n^p\subseteq\F_n^{p+1}\subseteq\cdots\subseteq\F_n^\infty.$$

Let $P$ be a path in $\F_n^p$. It can be written (uniquely)  
\begin{equation}\label{decomposition F_n^p}
    P=U^{i-1}QUD^i \text{ for some } i\in\{1,\ldots,p\} \text{ and } Q\in\F_{n-i}^p.
\end{equation}
Note that $P$ satisfies $\texttt{type}(P)=i$. This decomposition implies that   
$$|\F_n^p|=|\F_{n-1}^p|+|\F_{n-2}^p|+\ldots+|\F_{n-p}^p|,$$
which implies that $\F_n^p$ is enumerated by the generalized Fibonacci number $F_n^p$. Similarly, every path $P$ in $\F_n^\infty$ can be decomposed either $P=QUD$ or $P=UQD$ for some $Q\in\F_{n-1}^\infty$, which implies $|\F_n^\infty|=2^{n-1}$.

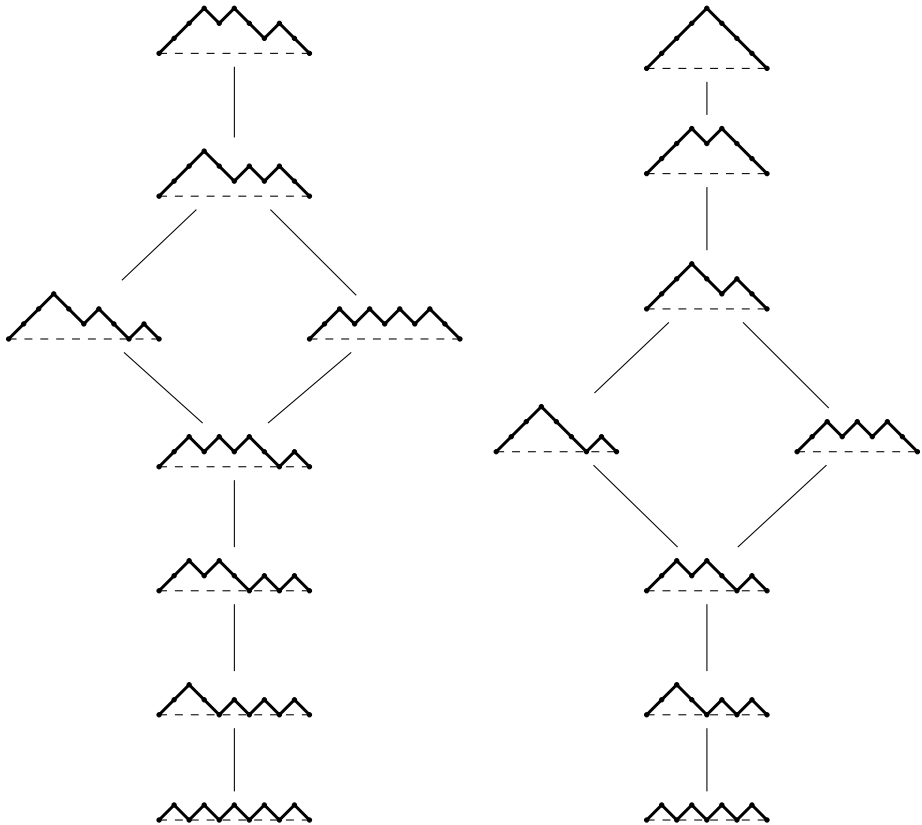
\begin{figure}[h]
\captionsetup{justification=centering,margin=1cm}
\begin{tikzpicture}
    \node (0) at (0,0) {\begin{tikzpicture}[scale=0.2]
        \draw[solid,line width=0.4mm] (0,0)-- ++\u -- ++\d -- ++\u -- ++\d -- ++\u -- ++\d -- ++\u -- ++\d -- ++\u -- ++\d;
         \filldraw (0,0)  circle (4pt);
        \filldraw (1,1)  circle (4pt);
        \filldraw (2,0)  circle (4pt);
        \filldraw (3,1)  circle (4pt);
        \filldraw (4,0)  circle (4pt);
        \filldraw (5,1)  circle (4pt);
        \filldraw (6,0)  circle (4pt);
        \filldraw (7,1)  circle (4pt);
        \filldraw (8,0)  circle (4pt);
        \filldraw (9,1)  circle (4pt);
        \filldraw (10,0)  circle (4pt);
        \draw[dashed] (0,0)--(10,0);
    \end{tikzpicture}};
    \node (1) at (0,1.5) {\begin{tikzpicture}[scale=0.2]
        \draw[solid,line width=0.4mm] (0,0)-- ++\u -- ++\u -- ++\d -- ++\d -- ++\u -- ++\d -- ++\u -- ++\d -- ++\u -- ++\d;
         \filldraw (0,0)  circle (4pt);
         \filldraw (1,1)  circle (4pt);
         \filldraw (2,2)  circle (4pt);
         \filldraw (3,1)  circle (4pt);
         \filldraw (4,0)  circle (4pt);
         \filldraw (5,1)  circle (4pt);
         \filldraw (6,0)  circle (4pt);
         \filldraw (7,1)  circle (4pt);
         \filldraw (8,0)  circle (4pt);
         \filldraw (9,1)  circle (4pt);
         \filldraw (10,0)  circle (4pt);
         \draw[dashed] (0,0)--(10,0);
    \end{tikzpicture}};
    \node (2) at (0,3.15) {\begin{tikzpicture}[scale=0.2]
        \draw[solid,line width=0.4mm] (0,0)-- ++\u -- ++\u -- ++\d -- ++\u -- ++\d -- ++\d -- ++\u -- ++\d -- ++\u -- ++\d;
         \filldraw (0,0)  circle (4pt);
         \filldraw (1,1)  circle (4pt);
         \filldraw (2,2)  circle (4pt);
         \filldraw (3,1)  circle (4pt);
         \filldraw (4,2)  circle (4pt);
         \filldraw (5,1)  circle (4pt);
         \filldraw (6,0)  circle (4pt);
         \filldraw (7,1)  circle (4pt);
         \filldraw (8,0)  circle (4pt);
         \filldraw (9,1)  circle (4pt);
         \filldraw (10,0)  circle (4pt);
         \draw[dashed] (0,0)--(10,0);
    \end{tikzpicture}};
    \node (3) at (0,4.8) {\begin{tikzpicture}[scale=0.2]
        \draw[solid,line width=0.4mm] (0,0)-- ++\u -- ++\u -- ++\d -- ++\u -- ++\d -- ++\u -- ++\d -- ++\d -- ++\u -- ++\d;
         \filldraw (0,0)  circle (4pt);
         \filldraw (1,1)  circle (4pt);
         \filldraw (2,2)  circle (4pt);
         \filldraw (3,1)  circle (4pt);
         \filldraw (4,2)  circle (4pt);
         \filldraw (5,1)  circle (4pt);
         \filldraw (6,2)  circle (4pt);
         \filldraw (7,1)  circle (4pt);
         \filldraw (8,0)  circle (4pt);
         \filldraw (9,1)  circle (4pt);
         \filldraw (10,0)  circle (4pt);
         \draw[dashed] (0,0)--(10,0);
    \end{tikzpicture}};
    \node (13) at (-2,6.6) {\begin{tikzpicture}[scale=0.2]
        \draw[solid,line width=0.4mm] (0,0)-- ++\u -- ++\u -- ++\u -- ++\d -- ++\d -- ++\u -- ++\d -- ++\d -- ++\u -- ++\d;
         \filldraw (0,0)  circle (4pt);
         \filldraw (1,1)  circle (4pt);
         \filldraw (2,2)  circle (4pt);
         \filldraw (3,3)  circle (4pt);
         \filldraw (4,2)  circle (4pt);
         \filldraw (5,1)  circle (4pt);
         \filldraw (6,2)  circle (4pt);
         \filldraw (7,1)  circle (4pt);
         \filldraw (8,0)  circle (4pt);
         \filldraw (9,1)  circle (4pt);
         \filldraw (10,0)  circle (4pt);
         \draw[dashed] (0,0)--(10,0);
    \end{tikzpicture}};
    \node (4) at (2,6.5) {\begin{tikzpicture}[scale=0.2]
        \draw[solid,line width=0.4mm] (0,0)-- ++\u -- ++\u -- ++\d -- ++\u -- ++\d -- ++\u -- ++\d -- ++\u -- ++\d -- ++\d;
         \filldraw (0,0)  circle (4pt);
         \filldraw (1,1)  circle (4pt);
         \filldraw (2,2)  circle (4pt);
         \filldraw (3,1)  circle (4pt);
         \filldraw (4,2)  circle (4pt);
         \filldraw (5,1)  circle (4pt);
         \filldraw (6,2)  circle (4pt);
         \filldraw (7,1)  circle (4pt);
         \filldraw (8,2)  circle (4pt);
         \filldraw (9,1)  circle (4pt);
         \filldraw (10,0)  circle (4pt);
         \draw[dashed] (0,0)--(10,0);
    \end{tikzpicture}};
    \node (14) at (0,8.5) {\begin{tikzpicture}[scale=0.2]
        \draw[solid,line width=0.4mm] (0,0)-- ++\u -- ++\u -- ++\u -- ++\d -- ++\d -- ++\u -- ++\d -- ++\u -- ++\d -- ++\d;
         \filldraw (0,0)  circle (4pt);
         \filldraw (1,1)  circle (4pt);
         \filldraw (2,2)  circle (4pt);
         \filldraw (3,3)  circle (4pt);
         \filldraw (4,2)  circle (4pt);
         \filldraw (5,1)  circle (4pt);
         \filldraw (6,2)  circle (4pt);
         \filldraw (7,1)  circle (4pt);
         \filldraw (8,2)  circle (4pt);
         \filldraw (9,1)  circle (4pt);
        \filldraw (10,0)  circle (4pt);
        \draw[dashed] (0,0)--(10,0);
    \end{tikzpicture}};
    \node (24) at (0,10.4) {\begin{tikzpicture}[scale=0.2]
        \draw[solid,line width=0.4mm] (0,0)-- ++\u -- ++\u -- ++\u -- ++\d -- ++\u -- ++\d -- ++\d -- ++\u -- ++\d -- ++\d;
         \filldraw (0,0)  circle (4pt);
         \filldraw (1,1)  circle (4pt);
         \filldraw (2,2)  circle (4pt);
         \filldraw (3,3)  circle (4pt);
         \filldraw (4,2)  circle (4pt);
         \filldraw (5,3)  circle (4pt);
         \filldraw (6,2)  circle (4pt);
         \filldraw (7,1)  circle (4pt);
         \filldraw (8,2)  circle (4pt);
         \filldraw (9,1)  circle (4pt);
         \filldraw (10,0)  circle (4pt);
         \draw[dashed] (0,0)--(10,0);
    \end{tikzpicture}};
    \draw (0) -- (1);
    \draw (1) -- (2);
    \draw (2) -- (3);
    \draw (3) -- (4);
    \draw (3) -- (13);
    \draw (4) -- (14);
    \draw (13) -- (14);
    \draw (14) -- (24);
\end{tikzpicture}
\begin{tikzpicture}
    \node (0) at (0,0) {\begin{tikzpicture}[scale=0.2]
        \draw[solid,line width=0.4mm] (0,0)-- ++\u -- ++\d -- ++\u -- ++\d -- ++\u -- ++\d -- ++\u -- ++\d;
        \filldraw (0,0)  circle (4pt);
        \filldraw (1,1)  circle (4pt);
        \filldraw (2,0)  circle (4pt);
       \filldraw (3,1)  circle (4pt);
        \filldraw (4,0)  circle (4pt);
       \filldraw (5,1)  circle (4pt);
        \filldraw (6,0)  circle (4pt);
       \filldraw (7,1)  circle (4pt);
        \filldraw (8,0)  circle (4pt);
        \draw[dashed] (0,0)--(8,0);
    \end{tikzpicture}};
    \node (1) at (0,1.5) {\begin{tikzpicture}[scale=0.2]
        \draw[solid,line width=0.4mm] (0,0)-- ++\u -- ++\u -- ++\d -- ++\d -- ++\u -- ++\d -- ++\u -- ++\d;
        \filldraw (0,0)  circle (4pt);
        \filldraw (1,1)  circle (4pt);
        \filldraw (2,2)  circle (4pt);
        \filldraw (3,1)  circle (4pt);
        \filldraw (4,0)  circle (4pt);
        \filldraw (5,1)  circle (4pt);
        \filldraw (6,0)  circle (4pt);
       \filldraw (7,1)  circle (4pt);
        \filldraw (8,0)  circle (4pt);
        \draw[dashed] (0,0)--(8,0);
    \end{tikzpicture}};
    \node (2) at (0,3.15) {\begin{tikzpicture}[scale=0.2]
        \draw[solid,line width=0.4mm] (0,0)-- ++\u -- ++\u -- ++\d -- ++\u -- ++\d -- ++\d -- ++\u -- ++\d;
        \filldraw (0,0)  circle (4pt);
        \filldraw (1,1)  circle (4pt);
        \filldraw (2,2)  circle (4pt);
        \filldraw (3,1)  circle (4pt);
       \filldraw (4,2)  circle (4pt);
        \filldraw (5,1)  circle (4pt);
        \filldraw (6,0)  circle (4pt);
        \filldraw (7,1)  circle (4pt);
        \filldraw (8,0)  circle (4pt);
        \draw[dashed] (0,0)--(8,0);
    \end{tikzpicture}};
    \node (12) at (-2,5.1) {\begin{tikzpicture}[scale=0.2]
        \draw[solid,line width=0.4mm] (0,0)-- ++\u -- ++\u -- ++\u -- ++\d -- ++\d -- ++\d -- ++\u -- ++\d;
        \filldraw (0,0)  circle (4pt);
        \filldraw (1,1)  circle (4pt);
        \filldraw (2,2)  circle (4pt);
        \filldraw (3,3)  circle (4pt);
        \filldraw (4,2)  circle (4pt);
        \filldraw (5,1)  circle (4pt);
        \filldraw (6,0)  circle (4pt);
        \filldraw (7,1)  circle (4pt);
        \filldraw (8,0)  circle (4pt);
        \draw[dashed] (0,0)--(8,0);
    \end{tikzpicture}};
    \node (3) at (2,5) {\begin{tikzpicture}[scale=0.2]
        \draw[solid,line width=0.4mm] (0,0)-- ++\u -- ++\u -- ++\d -- ++\u -- ++\d -- ++\u -- ++\d -- ++\d;
        \filldraw (0,0)  circle (4pt);
        \filldraw (1,1)  circle (4pt);
        \filldraw (2,2)  circle (4pt);
        \filldraw (3,1)  circle (4pt);
       \filldraw (4,2)  circle (4pt);
        \filldraw (5,1)  circle (4pt);
       \filldraw (6,2)  circle (4pt);
        \filldraw (7,1)  circle (4pt);
        \filldraw (8,0)  circle (4pt);
        \draw[dashed] (0,0)--(8,0);
    \end{tikzpicture}};
    \node (13) at (0,7) {\begin{tikzpicture}[scale=0.2]
        \draw[solid,line width=0.4mm] (0,0)-- ++\u -- ++\u -- ++\u -- ++\d -- ++\d -- ++\u -- ++\d -- ++\d;
        \filldraw (0,0)  circle (4pt);
        \filldraw (1,1)  circle (4pt);
        \filldraw (2,2)  circle (4pt);
        \filldraw (3,3)  circle (4pt);
        \filldraw (4,2)  circle (4pt);
        \filldraw (5,1)  circle (4pt);
        \filldraw (6,2)  circle (4pt);
        \filldraw (7,1)  circle (4pt);
        \filldraw (8,0)  circle (4pt);
        \draw[dashed] (0,0)--(8,0);
    \end{tikzpicture}};
    \node (23) at (0,8.8) {\begin{tikzpicture}[scale=0.2]
        \draw[solid,line width=0.4mm] (0,0)-- ++\u -- ++\u -- ++\u -- ++\d -- ++\u -- ++\d -- ++\d -- ++\d;
        \filldraw (0,0)  circle (4pt);
        \filldraw (1,1)  circle (4pt);
        \filldraw (2,2)  circle (4pt);
        \filldraw (3,3)  circle (4pt);
        \filldraw (4,2)  circle (4pt);
        \filldraw (5,3)  circle (4pt);
        \filldraw (6,2)  circle (4pt);
        \filldraw (7,1)  circle (4pt);
        \filldraw (8,0)  circle (4pt);
        \draw[dashed] (0,0)--(8,0);
    \end{tikzpicture}};
    \node (123) at (0,10.3) {\begin{tikzpicture}[scale=0.2]
        \draw[solid,line width=0.4mm] (0,0)-- ++\u -- ++\u -- ++\u -- ++\u -- ++\d -- ++\d -- ++\d -- ++\d;
        \filldraw (0,0)  circle (4pt);
        \filldraw (1,1)  circle (4pt);
        \filldraw (2,2)  circle (4pt);
        \filldraw (3,3)  circle (4pt);
        \filldraw (4,4)  circle (4pt);
        \filldraw (5,3)  circle (4pt);
        \filldraw (6,2)  circle (4pt);
        \filldraw (7,1)  circle (4pt);
        \filldraw (8,0)  circle (4pt);
        \draw[dashed] (0,0)--(8,0);
    \end{tikzpicture}};
    \draw (0) -- (1);
    \draw (1) -- (2);
    \draw (2) -- (3);
    \draw (2) -- (12);
    \draw (3) -- (13);
    \draw (12) -- (13);
    \draw (13) -- (23);
    \draw (23) -- (123);
\end{tikzpicture}   
    \caption{The Hasse diagrams of  $\FF_5^2$ (on the left) and $\FF_4^\infty$ (on the right).}
    \label{F_5^2 and F_4^infty}
\end{figure}

We equip $\F_n^p$ and $\F_n^\infty$ with the Stanley order $\leq$, and we denote by $\FF_n^p=(\F_n^p,\leq)$ and $\FF_n^\infty=(\F_n^\infty,\leq)$ the associated posets. See Figure~\ref{F_5^2 and F_4^infty} for an illustration of the Hasse diagrams of $\FF_5^2$ and $\FF_4^\infty$. We will write $P\lessdot Q$ when $Q$ \textit{covers} $P$, \textit{i.e.}, whenever $Q$ is obtained from $P$ by a transformation $DU\rightarrow UD$. We also say that $P$ is a \textit{lower cover} of $Q$, or $Q$ is an \textit{upper cover} of $P$. Then $\FF_n^p$ and $\FF_n^\infty$ are distributive lattices, as sublattices of the Stanley lattice. Indeed, let $P,Q\in\F_n^p$ (resp. $\F_n^\infty$), and let $P\wedge Q$ and $P\vee Q$ be respectively their meet (greatest lower bound) and their join (least upper bound) in the Stanley lattice. Then $P\wedge Q$ (resp. $P\vee Q$) is obtained by considering the lower (resp. upper) envelope of the two paths $P$ and $Q$, so they both clearly belong to $\F_n^p$ (resp. $\F_n^\infty$). Clearly, these lattices are ranked with the area below the path and above the $x$-axis.

We end this section by giving the classical concepts of partial order theory \cite{Grat} that we use in this study. A \textit{meet-irreducible} (resp. \textit{join-irreducible}) element is an element having exactly one upper (resp. lower) cover. An \textit{interval} $[P,Q]$ in a poset $\mathbb{P}$ is the set $\{R\in\mathbb{P}, P\leq R\leq Q\}$. The \textit{height} of $[P,Q]$ is the length of a maximal chain between $P$ and $Q$. An interval is said to be \textit{linear} when all its elements are pairwise comparable. An interval is \textit{boolean} when it is isomorphic to a boolean lattice. See \cite{BKN,BBKN,Bern,Bous,Boucha,Chap,chen1,chen2,Fan,Kre} for several studies on the enumeration of intervals into posets of Dyck paths.

\textbf{Outline of the paper.} In Section \ref{sec:covers and boolean} we enumerate the elements of $\FF_n^p$ ($p\geq2$) according to their number of upper-covers. We deduce from that the number of boolean intervals in $\FF_n^p$. We then use a discrete continuity argument to transfer those results to $\FF_n^\infty$. We also count meet-irreducible elements in $\FF_n^p$ which establishes a surprising link with the edges of the $(n,p)$-Tur\'an graph. In Section \ref{sec:linear intervals} we enumerate linear intervals, first in $\FF_n^p$ for $p\geq2$, and then in $\FF_n^\infty$, again using discrete continuity. In Section \ref{sec:intervals}, we prove that intervals are in one-to-one correspondence with bicolored Motzkin paths avoiding some patterns, which allows us to enumerate intervals for $p\in\{2,\infty\}$. Finally, we present in Section \ref{sec:bijections} bijections between the Dyck paths of $\F_n^p$ and other classical combinatorial objects. This gives new interpretations of the lattices $\FF_n^p$, in particular it can be seen as the well-known dominance order on some compositions.


\section{Coverings, irreducible elements, boolean intervals.}\label{sec:covers and boolean} 
In this section, we provide enumerative results for several characteristic elements (coverings, join- and meet-irreducible elements,  boolean intervals) in the lattices $\FF_n^p$, $p\geq 2$, and $\FF_n^\infty$.  


\subsection{In the lattices $\FF_n^p$, $p\geq 2$} Let $F_p(x,y)$ be the bivariate generating function where the coefficient of $x^ny^k$, $n,k\geq 0$, in its series expansion is the number of elements in  $\FF_n^p$ that have exactly $k$ upper covers. In this part we provide a closed form for $F_p(x,y)$, and we deduce the bivariate generating function $B_p(x,y)$ for the number of boolean intervals  in  $\FF_n^p$ with respect to $n$ and the interval height.

\begin{theorem} The generating function $F_p(x,y)$ is given by
  $$F_p(x,y)=\frac{(1-x)(1+(y-1)x^p)}{1-2x+x^{p + 1} - (y-1)(x^2-x^p+x^{p+1}-x^{p + 2})}.$$
  \label{covering}
\end{theorem}
\begin{proof} For the sake of consistency (only for this proof), even though it can look surprising at first, we consider that the empty path $\epsilon$ (the only element of $\F_0^p$) satisfies $\texttt{type}(\epsilon)=1$, and $UD$ (the only element of $\F_1^p$) satisfies $\texttt{type}(UD)=p$. 

For $1\leq i\leq p$, let $f_k^i(x)$ be the generating function for the number of elements of type $i$ in  $\FF_n^p$ having exactly $k$ upper covers, and we set $f_k(x)=\sum\limits_{i=1}^p f_k^i(x)$. We then have $$F_p(x,y)=\sum\limits_{k\geq 0}f_k(x)y^k.$$

With the above convention, the empty path has type 1 and is covered by $0$ element, so its contribution will appear in $f_0^1(x)$. Similarly, $UD$ has type $p$ and is covered by $0$ element, so its contribution will appear in  $f_0^p(x)$.

Let $P$ be a nonempty element in $\FF_n^p$ such that $\texttt{type}(P)=i$, and let $P=U^{i-1}QUD^i$ be its decomposition as described in (\ref{decomposition F_n^p}). Assume that $P$ has exactly $k$ upper covers.

\textbf{Case 1}: $P$ is of type $i\in[1,p-1]$. We distinguish two subcases ($a$) and ($b$).
\begin{enumerate}
    \item[($a$)]  $\texttt{type}(Q)=1$. With the above convention, either $Q$ is empty, or $Q$ ends with $UD$ and $Q\neq UD$. Then, $P$ and $Q$ have the same number of upper covers (in the case where $Q$ ends with $UD$, $Q\neq UD$, we cannot have a covering involving the last valley $DU$ of $P$ because the change $DU\rightarrow UD$ would create an occurrence $DUU$, see Figure~\ref{covering elements}).  So, the contribution of  these paths is $x^if_k^1(x)$.

\begin{figure}[h]
    \centering
    \captionsetup{justification=centering,margin=1cm}
    \begin{tikzpicture}[scale=0.4]
        \draw[solid,line width=0.4mm] (0,0)-- ++\u -- ++\u -- ++\u;
        \draw[solid,line width=0.4mm] (10,3)-- ++\u -- ++\d -- ++\d  -- ++\d -- ++\d;
        \draw[solid,line width=0.4mm,blue] (7,4)-- ++\d -- ++\u -- ++\d;
        \draw[solid,line width=0.4mm,red,dashed] (9,4)-- ++\u -- ++\d;
        \filldraw (0,0)  circle (4pt);
        \filldraw (1,1)  circle (4pt);
        \filldraw (2,2)  circle (4pt);
        \filldraw (3,3)  circle (4pt);
        \filldraw (10,3)  circle (4pt);
        \filldraw (11,4)  circle (4pt);
        \filldraw (12,3)  circle (4pt);
        \filldraw (13,2)  circle (4pt);
        \filldraw (14,1)  circle (4pt);
        \filldraw (15,0)  circle (4pt);
        \filldraw[blue] (8,3)  circle (4pt);
        \filldraw[blue] (9,4)  circle (4pt);
        \filldraw[blue] (7,4)  circle (4pt);
        \filldraw[red] (10,5)  circle (4pt);
        \node[blue] at (6.5,5.25) {$Q$};
        \draw[dashed] (0,0)--(15,0);
        \draw[dashed] (3,3)--(12,3);
        \draw[solid,line width=0.4mm,blue] plot [smooth, tension=1] coordinates {(3,3) (4,5) (5,3.2) (6,4.5) (7,4)};
    \end{tikzpicture}
\caption{Illustration of \textbf{Case~1}($a$) in the proof of Theorem~\ref{covering}. The last valley of $P$ cannot produce a covering because this would create an occurrence $DUU$.}
\label{covering elements}
\end{figure}
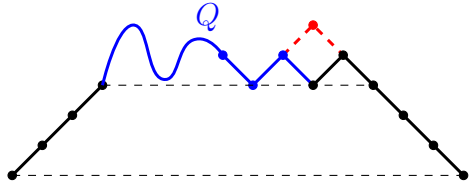

\item[($b$)] $\texttt{type}(Q)\in [2,p]$.  The last valley $DU$ of $P$ can produce a covering by replacing $DU$ with $UD$, and the $k-1$ remaining coverings are also coverings of $Q$. Note that with our convention, it is consistent with the case $Q=UD$. Then, the contribution of these paths is $x^i(f_{k-1}^2(x)+\ldots+f_{k-1}^p(x))$.
\end{enumerate}

\textbf{Case 2}: $P$ is of type $p$. Then $P$ and $Q$ have the same number of coverings, and the contribution of such paths is $x^p(f_k^1(x)+f_{k}^2(x)+\ldots+f_{k}^p(x))$.

Considering all these cases, we obtain the following system of equations for $k\geq 1$:
\begin{equation}\label{equations phi}
\left\{\begin{array}{ll}
    f_k^i(x) & =x^i(f_k^1(x)+f_{k-1}^2(x)+\ldots+f_{k-1}^p(x)) \text{ for } 1\leq i\leq p-1, \\
    f_k^p(x) & =x^p(f_k^1(x)+f_k^2(x)+\ldots+f_k^p(x)),
\end{array}\right.
\end{equation}
with the initial conditions 
$f_0^1(x)=1$, $f_0^p=\frac{x}{1-x}$, and $f_0^i(x)=0$ for $2\leq i\leq p-1$.
Indeed for every $n\geq 0$, the only path covered by 0 elements is the maximal one, and it has type $1$ only if $n=0$, otherwise it has type $p$ (remember our convention for $n=0$ and $n=1$). 

Now let $\varphi_k(x)=f_k^1(x)$ and $\theta_k(x)=f_k^2(x)+\ldots+f_k^p(x)$. From Equation (\ref{equations phi}), we get for $k\geq 1$
$$\left\{\begin{array}{ll}
    (1-x)\varphi_k(x) &= x\theta_{k-1}(x) \\
     (1-x^p)\theta_k(x) &= (x^2+\ldots+x^p)\varphi_k(x)+(x^2+\ldots+x^{p-1})\theta_{k-1}(x).  
\end{array}\right.$$
Observe that $\theta_0(x)=\frac{x}{1-x}$. Then we can solve for $\varphi_k$ and $\theta_k$ when $k\geq1$:

\begin{align*}   \varphi_k(x)&=\left(\frac{x^2-x^p+x^{p+1}-x^{p+2}}{(1-x)^2(1-x^p)}\right)^{k-1}\cdot\frac{x^2}{(1-x)^2},  \qquad \mbox{ and }\\
     \theta_k(x)&= \left(\frac{x^2-x^p+x^{p+1}-x^{p+2}}{(1-x)^2(1-x^p)}\right)^k\cdot\frac{x}{1-x}.
\end{align*}

We deduce that for all $k\geq1$, $$f_k(x)=\varphi_k(x)+\theta_k(x)=\left(\frac{x^2-x^p+x^{p+1}-x^{p+2}}{(1-x)^2(1-x^p)}\right)^{k-1}\cdot\frac{x^2-x^{p+1}}{(1-x)^3(1-x^p)}.$$
Finally, with $f_0(x)=\frac{1}{1-x}$, we obtain 
$$F_p(x,y)=\frac{1}{1-x}+\sum_{k=1}^{+\infty} f_k(x)y^k,$$
which yields the desired result after simplification.
\end{proof}

The first terms of the series expansions of $F_p(x,y)$ for $p=2$ and $p=3$ are respectively 
$$1 + x + (1+y)x^2 + (1+2y)x^3 + (1+4y)x^4 + (1+6y+y^2)x^5 + (1+9y+3y^2)x^6 +O(x^7)$$
and 
$$1 + x + (1+y)x^2 + (1+3y)x^3 + (1+5y+y^2)x^4 + (1+8y+4y^2)x^5 + (1+12y+10y^2+y^3)x^6 + O(x^7).$$
\begin{corollary}\label{nb coverings} The generating function for the number of coverings in $\FF_n^p$, $n\geq 0$,  is 
$$\partial_yF_p(x,y)\vert_{y=1}=\frac{\left(1-x \right)\left(x^2-x^{p+1}\right)\left(1-x^p\right)}{\left(1-2 x +x^{p +1}\right)^{2}}.$$
\end{corollary}
The first terms of the series expansion of $\partial_yF_2(x,y)\vert_{y=1}$ are 
$$x^2+2x^3+4x^4+8x^5+15x^6+28x^7+51x^8+92x^9+O(x^{10}),$$ and the sequence of coefficients corresponds to \seqnum{A029907} in \cite{oeis} where the $n$-th term $a_n$ satisfies 
$a_{n+1} = a_n + a_{n-1} + F_n$, where $a_0 =a_1= 0$ and $F_n$ is the $n$-th Fibonacci number (see Section \ref{sec:intro}). The sequences for $p\geq 3$ do not seem to appear in \cite{oeis}.

Moreover, the generating function $f_1(x)=\frac{x^2-x^{p+1}}{(1-x)^3(1-x^p)}$ (which is the coefficient of $y$ in the series expansion of $F_p(x,y)$) counts the number of meet-irreducible elements in $\FF_n^p$ (elements having exactly one upper cover). For $p=2$ the first terms of the series expansion are $$x^2+2x^3+4x^4+6x^5+9x^6+12x^7+16x^8+20x^9+25x^{10}+O(x^{11}).$$ The sequence of coefficients corresponds to \seqnum{A002620} in \cite{oeis}, where the $n$-th term is $b_2(n)=\lfloor\frac{n^2}{4}\rfloor$. More generally ($p\geq 2$), the following theorem provides a surprising link between the $n$-th coefficient $b_p(n)$ and a well-known parameter in extremal graph theory.

\begin{theorem}
For any $p\geq 2$, the number of meet-irreducible elements in $\FF_n^p$, that is the $n$-th coefficient of $f_1(x)$, is given by $$b_p(n)=\left\lfloor\frac{n^2(p-1)}{2p}\right\rfloor,$$  which also counts the number of edges in the $(n,p)$-Tur\'an graph (see \seqnum{A198787 } in \cite{oeis}, and \cite{Aig,bollobas,turan}).
\end{theorem}
\begin{proof}
    We present an analytic proof, showing that the generating function of $\left\lfloor\frac{n^2(p-1)}{2p}\right\rfloor$ is  $f_1(x)$. Considering the equality
    $$\left\lfloor\frac{n^2(p-1)}{2p}\right\rfloor=\left(1-\frac{1}{p}\right)\frac{n^2}{2}-\frac{(n\mod{p})(p-(n\mod p))}{2p},$$
    and observing that the generating functions of $(n\mod{p})$ and $(n\mod{p})^2$ satisfy     $$\sum_{n=0}^{+\infty} (n\mod{p})x^n=\left(\sum_{k=0}^{p-1}kx^k\right)\left(\sum_{k=0}^{+\infty} x^{pn}\right)=\frac{(p-1)x^{p+1}-px^p+x}{(1-x)^2(1-x^p)},$$
    and 
    \begin{align*}
        \sum_{n=0}^{+\infty} (n\mod{p})^2x^n&=\left(\sum_{k=0}^{p-1}k^2x^k\right)\left(\sum_{k=0}^{+\infty} x^{pn}\right)\\
        &=\frac{(2p^2 - 2p - 1)x^{p + 1} - (p - 1)^2x^{p + 2} - p^2x^p + x^2 + x}{(1-x)^3(1-x^p)},
    \end{align*}
    we can check that $\sum_{n=0}^{+\infty}\left\lfloor\frac{n^2(p-1)}{2p}\right\rfloor x^n=\frac{x^2-x^{p+1}}{(1-x)^3(1-x^p)}=f_1(x)$.
\end{proof}

From Theorem~\ref{covering}, we can easily deduce the following.
\begin{corollary}
    The generating function $B_p(x,y)$ for the number of boolean intervals in $\FF_n^p$,  with respect to the semilength $n\geq 0$, and the interval height is given by
    $$B_p(x,y)=\frac{(1-x)(1+yx^p)}{1-2x+x^{p + 1} -y(x^2-x^p+x^{p+1}-x^{p + 2})}.$$
    \label{corbool}
\end{corollary}
\begin{proof} 

As a direct consequence of the distributivity of $\FF_n^p$, we have that $B_p(x,y)=F_p(x,1+y)$, see for instance \cite[Exercise 3.19]{enum comb}.
\end{proof}

The first terms of the series expansion of $B_2(x,y)$  are  
$$ 1+x+(2+y)x^2+(3+2y)x^3+(5+4y)x^4+(8+8y+y^2)x^5+(13+15y+3y^2)x^6+O(x^7).$$
For $p=2$, we find $B_2(x,1)=\frac{1+x^2}{1-x-x^2-x^3}$, so the boolean intervals of $\FF_n^2$ are enumerated by the Tribonacci numbers (\seqnum{A000213} in \cite{oeis}). For $p=3$, the number sequence of boolean intervals in $\FF_n^3$ also appears in \cite{oeis} (see \seqnum{A193641}). It seems to be the only two values of $p\geq 2$ such that the sequence appears in \cite{oeis}.

\begin{corollary}\label{symmetry degree}
    The coefficient of $x^ny^k$ in $F_p(x,y)$ is also the number of elements in $\FF_n^p$ that have exactly $k$ lower covers.
\end{corollary}
\begin{proof}
    Let $G_p(x,y)$ be the generating function whose coefficient of $x^ny^k$ is the number of elements in $\FF_n^p$ that have exactly $k$ lower covers. Again by the distributivity of $\FF_n^p$ \cite[Exercise 3.19]{enum comb}, we have $B_p(x,y)=G_p(x,1+y)=F_p(x,1+y)$, so $G_p(x,y)=F_p(x,y)$.
\end{proof}

\begin{remark}
For $(k,\ell)\in\N^2$, we say that an element $P\in\FF_n^p$ has \textit{degree} $(k,\ell)$ if $P$ has $k$ upper covers and $\ell$ lower covers. Given Corollary \ref{symmetry degree}, we could ask if there are as many elements with degree $(k,\ell)$ as elements with degree $(\ell,k)$ in $\FF_n^p$. Actually this is not the case for infinitely many values of $(n,p)$. Indeed, for $n>p$ and $(n\mod{p})\not\in\{0,1\}$, the maximal element has degree $(0,2)$, but the minimal element has degree $(1,0)$. However those two elements are the only ones whose degree has a $0$-coordinate. We will see later with Theorem~\ref{symmetry F_n^infty} that the distribution of the degree is always symmetric in $\FF_n^\infty$.
\end{remark}


\subsection{In the lattice $\FF_n^\infty$}
 Let $F_\infty(x,y)$ be the bivariate generating function where the coefficient of $x^ny^k$, $n,k\geq 0$, in its series expansion is the number of elements in  $\FF_n^\infty$ that have exactly $k$ upper covers, and let $B_\infty(x,y)$ be the bivariate generating function for the number of boolean intervals  in  $\FF_n^\infty$ with respect to $n$ and the interval height.
 
It is possible to make a similar study for $\FF_n^\infty$ as in the proofs of  Theorem~\ref{covering}, and Corollary~\ref{corbool}. However, we can also use a discrete continuity argument. Indeed, for any $P\in\F_n^\infty$, there exists $p\geq 2$ sufficiently large such that $P\in\F_n^p$. In fact, for any $n\geq 0$, we have $\F_n^\infty=\F_n^n$. Thus, if $\text{val}_x(\cdot)$ denotes the valuation of a power series in the variable $x$, then $\text{val}_x(F_\infty(x,y)-F_p(x,y))\geq p$ for all $p\geq 2$. Thus, $F_\infty(x,y)=\lim_{p\to\infty}F_p(x,y)$ (relatively to the metric $\text{val}_x(\cdot)$), and $B_\infty(x,y)=\lim_{p\to\infty}B_p(x,y)$.
Therefore, we have the following.

\begin{corollary}\label{boolean F infty} The generating functions $F_\infty(x,y)$ and $B_\infty(x,y)$ are given by
   $$F_\infty(x,y)=\frac{1-x}{1-2x+(1-y)x^2},\qquad \mbox{ and }\qquad
    B_\infty(x,y)=\frac{1-x}{1-2x-x^2y}.$$
    The number of coverings is  
    $$[x^n]\partial_yF_\infty(x,y)\vert_{y=1}=n\cdot2^{n-3} \ \text{ if } n\geq2, \text{ and } 0 \text{ otherwise,}$$ which corresponds to \seqnum{A001792} in \cite{oeis}.
    
    The number of boolean intervals in $\FF_n^\infty$ is $$[x^n]B_\infty(x,1)=\sum_{k=0}^{\lfloor n/2\rfloor}\binom{n}{n-2k}2^k, $$
    which corresponds to \seqnum{A078057} in \cite{oeis}.
    
    The number of meet-irreducible elements (resp. join-irreducible) is $\frac{n(n-1)}{2}$.
\end{corollary}

The first terms of the series expansions of $F_\infty(x,y)$ and $B_\infty(x,y)$ are respectively $$1 + x + (1 + y)x^2 + (1 + 3y)x^3 + (1+6y+y^2)x^4 + (1+10y+5y^2)x^5 + (1+15y+15y^2+y^3)x^6+O(x^7),$$
and
$$1 + x + (2+y)x^2 + (4+3y)x^3 + (8+8y+y^2)x^4 + (16+20y+5y^2)x^5 + (32+48y+18y^2+y^3)x^6+O(x^7).$$

\begin{theorem}
    For all $n\geq0,$ let $B_n$ be a random variable following the uniform distribution over the boolean intervals of $\FF_n^\infty$, and let $X_n$ be the height of $B_n$. Then $\frac{4X_n-(2-\sqrt{2})n}{\sqrt{n}\sqrt[4]{2}}$ converges in law to a standard normal distribution.
\end{theorem}
\begin{proof}
    We use singularity analysis, see \cite{analytic combinatorics}. Let $y$ belong to a neighbourhood of $1$. The main singularity of $B_\infty(x,y)$ is then $\frac{\sqrt{1+y}-1}{y}$, and near this singularity, we have the approximation
    $$B_\infty(x,y)\sim \frac{1}{2}\left(1-\frac{yx}{\sqrt{1+y}-1}\right)^{-1}.$$
    Since $[x^n]B_\infty(x,1)\sim\frac{(\sqrt{2}-1)^{-n}}{2}$, we have
    $$\frac{[x^n]B_\infty(x,y)}{[x^n]B_\infty(x,1)}\sim\left(\frac{y(\sqrt{2}-1)}{\sqrt{1+y}-1}\right)^n,$$
    then $X_n$ can be approximated by a sum of independent random variables with generating function $y\mapsto\frac{y(\sqrt{2}-1)}{\sqrt{1+y}-1}$. Since the last distribution has expected value $(2-\sqrt{2})/4$ and standard deviation $\sqrt[4]{2}/4$, by the central limit theorem, we have the convergence in law of the standardized version of $X_n$ to a standard normal distribution.
\end{proof}

\begin{definition}
    The Möbius function $\mu$ on a poset $\mathbb{P}$ is defined \cite{blass sagan,BFT} recursively by
    $$\mu(P,Q)=\left\{\begin{array}{cc}
        0 & \text{ if } P\not\leq Q,  \\
        1 & \text{ if } P=Q, \text{ and }\\
        -\sum_{P\leq R < Q} \mu(P,R) & \text{ for all } P<Q.
    \end{array}\right.$$
\end{definition}

Since $\FF_n^p$ is a finite distributive lattice, it is well-known (see for instance \cite{enum comb}) that, for any $P,Q\in\FF_n^p$, the Möbius function $\mu(P,Q)$ is equal to 0 if the interval $[P,Q]$ is not boolean, and otherwise $\mu(P,Q)=(-1)^h$, where $h$ is the height of $[P,Q]$. Again by the distributivity, if $[P,Q]$ is boolean, then $Q$ is the join of the upper covers of $P$. This implies that there exist factors $\alpha_1,\ldots,\alpha_{h+1}$ such that $P=\alpha_1DU\alpha_2DU\alpha_3\ldots \alpha_kDU\alpha_{h+1}$ and $Q=\alpha_1UD\alpha_2UD\alpha_3\ldots \alpha_kUD\alpha_{h+1}$. In terms of the rank  function $\rho$ of the lattice (namely the area between the path and the $x$-axis, as mentioned in Section \ref{sec:intro}), $h=\rho(Q)-\rho(P)$. See also Remark \ref{rank} for an alternative definition of the rank.


\section{Linear intervals}\label{sec:linear intervals}

In this section, we focus on the enumeration of the linear intervals in $\FF_n^p$, $p\geq 2$, and $\FF_n^\infty$. We first give three lemmas that characterize the structure of linear intervals. Next, we deduce the generating function for the number of linear intervals in $\FF_n^p$ with respect to $n$ and the interval height. Since degenerate cases arise only for $p=2$, the enumerations for $p=2$ and $p\geq3$ are handled separately in two distinct subsections. The generating function for the case $p=\infty$ is obtained using a discrete continuity argument.

\begin{lemma}\label{structure lin int F_n^p j=i}
    Let $P,Q\in\FF_n^p$, $p\geq 2$, and $1\leq i\leq p$, such that $i=\texttt{type}(P)=\texttt{type}(Q)$. Then we have  $[P,Q]=[U^{i-1}P'UD^i,U^{i-1}Q'UD^i]$ and $[P,Q]$ is a linear interval if and only if $[P',Q']$ is a linear interval in $\F_{n-i}^p$.
\end{lemma}

\begin{figure}[h]
    \centering
    \captionsetup{justification=centering,margin=1cm}
    \begin{tikzpicture}[scale=0.5]
        \draw[solid,line width=0.4mm,blue] (0,0)-- ++\u -- ++\u -- ++\u;
        \draw[solid,line width=0.4mm,blue] (10,3)-- ++\u -- ++\d -- ++\d  -- ++\d-- ++\d;
        \draw[solid,line width=0.4mm,blue] plot [smooth, tension=1] coordinates {(3,3) (5,5) (7,3.2) (8,4.5) (10,3)};
        \filldraw[blue] (0,0)  circle (4pt);
        \filldraw[blue] (1,1)  circle (4pt);
        \filldraw[blue] (2,2)  circle (4pt);
        \filldraw[blue] (3,3)  circle (4pt);
        \filldraw[blue] (10,3)  circle (4pt);
        \filldraw[blue] (11,4)  circle (4pt);
        \filldraw[blue] (12,3)  circle (4pt);
        \filldraw[blue] (13,2)  circle (4pt);
        \filldraw[blue] (14,1)  circle (4pt);
        \filldraw[blue] (15,0)  circle (4pt);
        \draw[blue] (8.3,3.6) node {$P'$};
         \draw (14,2.4) node {$D^i$};
        \draw[solid,line width=0.4mm,red] (0,0.3)-- ++\u -- ++\u -- ++\u;
        \draw[solid,line width=0.4mm,red] (10,3.3)-- ++\u -- ++\d -- ++\d  -- ++\d -- ++\d;
        \draw[solid,line width=0.4mm,red] plot [smooth, tension=1] coordinates {(3,3.3) (5,5.5) (7,4.3) (8,5.3) (10,3.3)};
        \draw[dashed] (10,3.15)--(3,3.15);
        \draw[dashed] (0,0)--(15,0);
        \filldraw[red] (0,0.3)  circle (4pt);
        \filldraw[red] (1,1.3)  circle (4pt);
        \filldraw[red] (2,2.3)  circle (4pt);
        \filldraw[red] (3,3.3)  circle (4pt);
        \filldraw[red] (10,3.3)  circle (4pt);
        \filldraw[red] (11,4.3)  circle (4pt);
        \filldraw[red] (12,3.3)  circle (4pt);
        \filldraw[red] (13,2.3)  circle (4pt);
        \filldraw[red] (14,1.3)  circle (4pt);
        \filldraw[red] (15,0.3)  circle (4pt);
        \draw[red] (8.7,5.6) node {$Q'$};
        \draw (7.5,1) node {Structure (1)};
    \end{tikzpicture}
    \caption{The structure of linear intervals $[P,Q]$ in $\FF_n^p$ when $\texttt{type}(P)=\texttt{type}(Q)$. See Lemma~\ref{structure lin int F_n^p j=i}}.
    \label{fig:structure lin int j=i}
\end{figure}

\begin{proof}
Since $P$ and $Q$ have the same type $i$, we can decompose $P=U^{i-1}P'UD^i$ and $Q=U^{i-1}Q'UD^i$, with $P',Q'\in\F_{n-i}^p$. Clearly, the two intervals $[P',Q']$ and $[P,Q]$ are isomorphic as posets, which complete the proof.  See Figure~\ref{fig:structure lin int j=i}.
\end{proof}

\begin{lemma}\label{structure lin int F_n^p j-i > 2}
    Let $P,Q\in\FF_n^p$, $p\ge 3$, such that $\texttt{type}(Q)\geq \texttt{type}(P)+2$. If $i:=\texttt{type}(P)$ and $j:=\texttt{type}(Q)$, then $[P,Q]$ is a linear interval if and only if
        \begin{enumerate}[label=(\alph*)]
            \item $[P,Q]=[U^{n-3}(UD)^3D^{n-3},U^{n}D^{n}]$, $3\leq n\leq p$, or
             \item $[P,Q]=[U^{j-1}RD^{j-i}UD^i, U^{j-1}RUD^j]$, with $R$ empty or $R\in\F_{n-j}^p$ and $\texttt{type}(R)\in[1,p-(j-i)]$.
        \end{enumerate}
\end{lemma}

\begin{figure}[h]
    \centering
    \captionsetup{justification=centering,margin=1cm}
    \begin{tikzpicture}[scale=0.5]
        \draw[solid,line width=0.4mm,blue] (0,0)-- ++\u -- ++\u -- ++\u-- ++\u -- ++\u -- ++\d -- ++\u -- ++\d -- ++\u -- ++\d -- ++\d-- ++\d -- ++\d  -- ++\d;
        \filldraw[blue] (0,0)  circle (4pt);
        \filldraw[blue] (1,1)  circle (4pt);
        \filldraw[blue] (2,2)  circle (4pt);
        \filldraw[blue] (3,3)  circle (4pt);
        \filldraw[blue] (4,4)  circle (4pt);
        \filldraw[blue] (5,5)  circle (4pt);
        \filldraw[blue] (6,4)  circle (4pt);
        \filldraw[blue] (7,5)  circle (4pt);
        \filldraw[blue] (8,4)  circle (4pt);
        \filldraw[blue] (9,5)  circle (4pt);
        \filldraw[blue] (10,4)  circle (4pt);
        \filldraw[blue] (11,3)  circle (4pt);
        \filldraw[blue] (12,2)  circle (4pt);
        \filldraw[blue] (13,1)  circle (4pt);
        \filldraw[blue] (14,0)  circle (4pt);
        \draw[dashed] (0,0)--(14,0);
        \draw[dashed] (4,4.1)--(10,4.1);
        \draw[solid,line width=0.4mm,red] (0,0.3)-- ++\u -- ++\u -- ++\u-- ++\u -- ++\u -- ++\u -- ++\u -- ++\d -- ++\d  -- ++\d -- ++\d-- ++\d -- ++\d  -- ++\d;
        \filldraw[red] (0,0.3)  circle (4pt);
        \filldraw[red] (1,1.3)  circle (4pt);
        \filldraw[red] (2,2.3)  circle (4pt);
        \filldraw[red] (3,3.3)  circle (4pt);
        \filldraw[red] (4,4.3)  circle (4pt);
        \filldraw[red] (5,5.3)  circle (4pt);
        \filldraw[red] (6,6.3)  circle (4pt);
        \filldraw[red] (7,7.3)  circle (4pt);
        \filldraw[red] (8,6.3)  circle (4pt);
        \filldraw[red] (9,5.3)  circle (4pt);
        \filldraw[red] (10,4.3)  circle (4pt);
        \filldraw[red] (11,3.3)  circle (4pt);
        \filldraw[red] (12,2.3)  circle (4pt);
        \filldraw[red] (13,1.3)  circle (4pt);
        \filldraw[red] (14,0.3)  circle (4pt);
        \draw (7,1) node {Structure (2$a$)};
    \end{tikzpicture}\qquad
    \begin{tikzpicture}[scale=0.5]
        \draw[solid,line width=0.4mm,blue] (0,0)-- ++\u -- ++\u -- ++\u;
        \draw[solid,line width=0.4mm,blue] (10,3)-- ++\d -- ++\d -- ++\u  -- ++\d -- ++\d;
        \draw[solid,line width=0.6mm] plot [smooth, tension=1] coordinates {(3,3) (5,5.5) (7,3.7) (8,5) (10,3)};
        \filldraw[blue] (0,0)  circle (4pt);
        \filldraw[blue] (1,1)  circle (4pt);
        \filldraw[blue] (2,2)  circle (4pt);
        \filldraw[blue] (3,3)  circle (4pt);
        \filldraw[blue] (10,3)  circle (4pt);
        \filldraw[blue] (11,2)  circle (4pt);
        \filldraw[blue] (12,1)  circle (4pt);
        \filldraw[blue] (13,2)  circle (4pt);
        \filldraw[blue] (14,1)  circle (4pt);
        \filldraw[blue] (15,0)  circle (4pt);
        \draw[solid,line width=0.4mm,red] (0,0.3)-- ++\u -- ++\u -- ++\u;
        \draw[solid,line width=0.4mm,red] (10,3.3)-- ++\u -- ++\d -- ++\d  -- ++\d -- ++\d;
         \draw[dashed] (0,0)--(15,0);
        \draw[dashed] (3,3.15)--(10,3.15);
        
        \filldraw[red] (0,0.3)  circle (4pt);
        \filldraw[red] (1,1.3)  circle (4pt);
        \filldraw[red] (2,2.3)  circle (4pt);
        \filldraw[red] (3,3.3)  circle (4pt);
        \filldraw[red] (10,3.3)  circle (4pt);
        \filldraw[red] (11,4.3)  circle (4pt);
        \filldraw[red] (12,3.3)  circle (4pt);
        \filldraw[red] (13,2.3)  circle (4pt);
        \filldraw[red] (14,1.3)  circle (4pt);
        \filldraw[red] (15,0.3)  circle (4pt);
        \draw (8.7,5.6) node {$R$};
        \draw (7.5,1) node {Structure ($2b$)};
    \end{tikzpicture}
    \caption{The structure of linear intervals $[P,Q]$ in $\FF_n^p$ when  $\texttt{type}(Q)\geq\texttt{type}(P)+2$. See  Lemma~\ref{structure lin int F_n^p j-i > 2}}
    \label{fig:strucure lin int j-i > 2}
\end{figure}
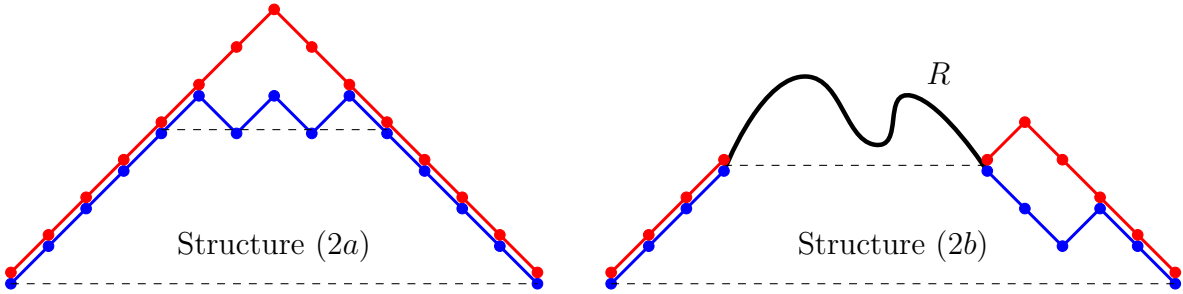

\begin{proof}
It is direct to check that intervals of the form \textit{(a)} or \textit{(b)} are linear (see Figure~\ref{fig:strucure lin int j-i > 2}), so we focus on the converse. Let $[P,Q]$ be a linear interval such that $i:=\texttt{type}(P)$ and $j:=\texttt{type}(Q)$ with $j\geq i+2$. We can write $P=U^{i-1}P'D^{i-1}$, where $P'\in\F_{n-i+1}^p$ is nonempty and $\texttt{type}(P')=1$. Similarly, $Q=U^{i-1}Q'D^{i-1}$, where $Q'\in\F_{n-i+1}$ is nonempty and $\texttt{type}(Q')=j-i+1$. Note that the intervals $[P,Q]$ and $[P',Q']$ are isomorphic as posets. So, let us distinguish four cases depending on the form of $P'$.

\begin{enumerate}[label=(\roman*)]
    \item Suppose that $P'=(UD)^{n-i+1}$. Since $j-i\geq2$, we have necessarily $n-i+1\geq3$. By a simple observation, if $n-i+1=3$, then $[P,Q]$ has the form \textit{(a)}. If $n-i+1\geq 4$, then $[P',Q']$ cannot be a linear interval.  Indeed, $U(UD)^{n-i}D$ and $U^3D^2UD^2(UD)^{n-i-3}$ both belong to $[P',Q']$ (they are obviously greater than $P'$, and they are lower than $Q'$ because $j-i\geq2$), and these two paths are not comparable, which condemns the linearity of $[P',Q']$, and so the one of $[P,Q]$. 
    \item Suppose that there exist $k\geq3$, $k'\geq2$ and a prefix $S$ such that $P'=SD^k(UD)^{k'}$, then $[P',Q']$ cannot be a linear interval. Indeed, the paths $SD^{k-2}UD^3(UD)^{k'-1}$ and $SD^{k-1}(UD)^2D(UD)^{k'-2}$ both belong to $[P',Q']$ and they are not comparable.
    \item Suppose that there exist $k'\geq2$ and a prefix $S$ such that $P'=SDUD^2(UD)^{k'}$,  then $[P',Q']$ cannot be a linear interval. Indeed, $SUD^3(UD)^{k'}$ and $SD(UD)^2D(UD)^{k'-1}$ both belong to $[P',Q']$ but they are not comparable.
    \item Suppose that $P'=U^2D^2(UD)^{n-i-1}$,  then $[P',Q']$ cannot be a linear interval. Indeed, $U^3D^3(UD)^{n-i-2}$ and $U(UD)^3D(UD)^{n-i-2}$ both belong to $[P',Q']$ but they are not comparable.
\end{enumerate}

We have proved that either $[P,Q]$ has the form \textit{(a)}, or  $P'$ has the  form $P'=SUD^kUD$ with $2\leq k\leq j-i$ (\textit{i.e.} $k'=1$ with this above notation).

If $k<j-i$, $[P',Q']$ cannot be a linear interval. Indeed , in this case we can write $P'=TDUD^kUD$, and the paths $TDUD^{k+1}UD$ and $TDUD^{k-1}UD^2$ both belong to $[P',Q']$, and they are not comparable. 

If $k=j-i$, then $P$ can be decomposed $P=U^{j-1}P'D^{j-i}UD^i$, with $P'\in\F_{n-j}^p$. Since $Q$ has type $j$ we also have $Q=U^{j-1}Q'D^{j-i}UD^i$, with $Q'\in\F_{n-j}^p$. Now if $P'\ne Q'$, $[P',Q']$ is a non trivial interval in $\F_{n-j}^p$, so there exists $R'\in [P',Q']$ covering $P'$. Then $U^{j-1}R'D^{j-i}UD^i$ and $U^{j-1}P'D^{j-i-1}UD^{i+1}$ both belong to $[P,Q]$, but they are not comparable, contradicting once again the linearity. We conclude that $P'=Q'$, so $[P,Q]$ has the form \textit{(b)}.
\end{proof}

\begin{lemma}\label{structure lin int F_n^p j-i=1}
    Let $P,Q\in\FF_n^p$, $p\geq 2$, such that $\texttt{type}(Q)= \texttt{type}(P)+1$. If $i:=\texttt{type}(P)$ and $j:=\texttt{type}(Q)$, then $[P,Q]$ is a linear interval if and only if 
        \begin{enumerate}[label=(\alph*)]
            \item $[P,Q]=[U^iR(DU)^kD^i,U^iR(UD)^kD^i]$ with $k\geq1$, $R$ empty, or $R\in\F_{n-k-i}^p$ and $\texttt{type}(R)\in[1,p-1]$, or
            \item $[P,Q]=[U^{k(p-1)+i}RD(UD^p)^kUD^i,U^{k(p-1)+i}R(UD^p)^kUD^{i+1}]$, with $k\geq1$, $R$ empty or $R\in\F_{n-kp-i-1}^p$, and and $\texttt{type}(R)\in[1,p-1]$, or
            \item $[P,Q]=[U^2D^2(UD)^2,U^3D^2UD^2]$, or $[P,Q]=[(UD)^4,U^3D^2UD^2]$, with $p=2$.
        \end{enumerate}
\end{lemma}

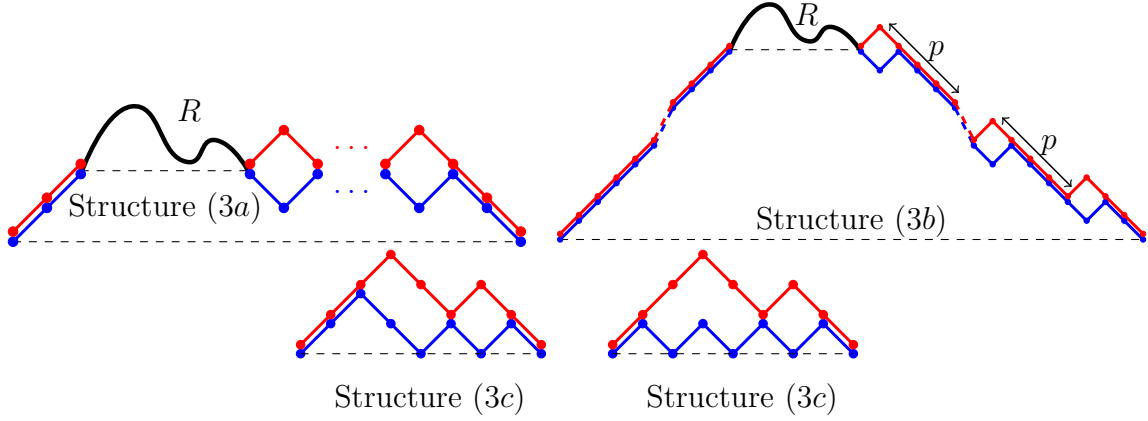
\begin{figure}[h]
    \centering
    \captionsetup{justification=centering,margin=1cm}
    \begin{tikzpicture}[scale=0.45]
        \draw[solid,line width=0.4mm,blue] (0,0)-- ++\u -- ++\u;
        \draw[solid,line width=0.4mm,blue] (7,2)-- ++\d -- ++\u;
        \draw[solid,line width=0.4mm,blue] (11,2)-- ++\d  -- ++\u -- ++\d -- ++\d;
        \draw[solid,line width=0.6mm] plot [smooth, tension=1] coordinates {(2,2) (3.5,4) (5,2.4) (6,3) (7,2)};
        \draw[dashed] (0,0)--(15,0);
        \draw[dashed] (2,2.1)--(7,2.1);
        \filldraw[blue] (0,0)  circle (4pt);
        \filldraw[blue] (1,1)  circle (4pt);
        \filldraw[blue] (2,2)  circle (4pt);
        \filldraw[blue] (7,2)  circle (4pt);
        \filldraw[blue] (8,1)  circle (4pt);
        \filldraw[blue] (9,2)  circle (4pt);
        \filldraw[blue] (11,2)  circle (4pt);
        \filldraw[blue] (12,1)  circle (4pt);
        \filldraw[blue] (13,2)  circle (4pt);
        \filldraw[blue] (14,1)  circle (4pt);
        \filldraw[blue] (15,0)  circle (4pt);
        \draw[blue] (10,1.5) node {$\ldots$};
        \draw[solid,line width=0.4mm,red] (0,0.3)-- ++\u -- ++\u;
        \draw[solid,line width=0.4mm,red] (7,2.3)-- ++\u -- ++\d;
        \draw[solid,line width=0.4mm,red] (11,2.3) -- ++\u  -- ++\d -- ++\d -- ++\d;
        \filldraw[red] (0,0.3)  circle (4pt);
        \filldraw[red] (1,1.3)  circle (4pt);
        \filldraw[red] (2,2.3)  circle (4pt);
        \filldraw[red] (7,2.3)  circle (4pt);
        \filldraw[red] (8,3.3)  circle (4pt);
        \filldraw[red] (9,2.3)  circle (4pt);
        \filldraw[red] (11,2.3)  circle (4pt);
        \filldraw[red] (12,3.3)  circle (4pt);
        \filldraw[red] (13,2.3)  circle (4pt);
        \filldraw[red] (14,1.3)  circle (4pt);
        \filldraw[red] (15,0.3)  circle (4pt);
        \draw[red] (10,2.8) node {$\ldots$};
        \draw (5.2,3.9) node {$R$};
        \draw (4.5,1) node {Structure ($3a$)};
    \end{tikzpicture}\quad
    \begin{tikzpicture}[scale=0.25]
        \draw[solid,line width=0.6mm] plot [smooth, tension=1] coordinates {(9,10) (11,12.5) (13,10.6) (14.5,11.3) (16,10)};
        \draw[dashed] (9,10.1)--(16,10.1);
        \draw[dashed] (0,0)--(31,0);
        \draw[solid,line width=0.4mm,blue] (0,0)-- ++\u -- ++\u -- ++\u-- ++\u -- ++\u;
        \draw[solid,line width=0.4mm,blue] (6,7)-- ++\u -- ++\u -- ++\u;
        \draw[solid,line width=0.4mm,blue] (16,10)-- ++\d -- ++\u -- ++\d -- ++\d -- ++\d;
        \draw[solid,line width=0.4mm,blue] (22,5)-- ++\d -- ++\u -- ++\d -- ++\d -- ++\d -- ++\d -- ++\u -- ++\d -- ++\d;
        \draw[dashed,line width=0.4mm,blue] (5,5)-- (6,7);
        \draw[dashed,line width=0.4mm,blue] (21,7)-- (22,5);
        \filldraw[blue] (0,0)  circle (4pt);
        \filldraw[blue] (1,1)  circle (4pt);
        \filldraw[blue] (2,2)  circle (4pt);
        \filldraw[blue] (3,3)  circle (4pt);
        \filldraw[blue] (4,4)  circle (4pt);
        \filldraw[blue] (5,5)  circle (4pt);
        \filldraw[blue] (6,7)  circle (4pt);
        \filldraw[blue] (7,8)  circle (4pt);
        \filldraw[blue] (8,9)  circle (4pt);
        \filldraw[blue] (9,10)  circle (4pt);
        \filldraw[blue] (16,10)  circle (4pt);
        \filldraw[blue] (17,9)  circle (4pt);
        \filldraw[blue] (18,10)  circle (4pt);
        \filldraw[blue] (19,9)  circle (4pt);
        \filldraw[blue] (20,8)  circle (4pt);
        \filldraw[blue] (21,7)  circle (4pt);
        \filldraw[blue] (22,5)  circle (4pt);
        \filldraw[blue] (23,4)  circle (4pt);
        \filldraw[blue] (24,5)  circle (4pt);
        \filldraw[blue] (25,4)  circle (4pt);
        \filldraw[blue] (26,3)  circle (4pt);
        \filldraw[blue] (27,2)  circle (4pt);
        \filldraw[blue] (28,1)  circle (4pt);
        \filldraw[blue] (29,2)  circle (4pt);
        \filldraw[blue] (30,1)  circle (4pt);
        \filldraw[blue] (31,0)  circle (4pt);
        \draw[solid,line width=0.4mm,red] (0,0.3)-- ++\u -- ++\u -- ++\u-- ++\u -- ++\u;
        \draw[solid,line width=0.4mm,red] (6,7.3)-- ++\u -- ++\u -- ++\u;
        \draw[solid,line width=0.4mm,red] (16,10.3)-- ++\u -- ++\d -- ++\d -- ++\d -- ++\d;
        \draw[solid,line width=0.4mm,red] (22,5.3)-- ++\u -- ++\d -- ++\d -- ++\d -- ++\d-- ++\u -- ++\d -- ++\d -- ++\d;
        \draw[dashed,line width=0.4mm,red] (5,5.3)-- (6,7.3);
        \draw[dashed,line width=0.4mm,red] (21,7.3)-- (22,5.3);
        \filldraw[red] (0,0.3)  circle (4pt);
        \filldraw[red] (1,1.3)  circle (4pt);
        \filldraw[red] (2,2.3)  circle (4pt);
        \filldraw[red] (3,3.3)  circle (4pt);
        \filldraw[red] (4,4.3)  circle (4pt);
        \filldraw[red] (5,5.3)  circle (4pt);
        \filldraw[red] (6,7.3)  circle (4pt);
        \filldraw[red] (7,8.3)  circle (4pt);
        \filldraw[red] (8,9.3)  circle (4pt);
        \filldraw[red] (9,10.3)  circle (4pt);
        \filldraw[red] (16,10.3)  circle (4pt);
        \filldraw[red] (17,11.3)  circle (4pt);
        \filldraw[red] (18,10.3)  circle (4pt);
        \filldraw[red] (19,9.3)  circle (4pt);
        \filldraw[red] (20,8.3)  circle (4pt);
        \filldraw[red] (21,7.3)  circle (4pt);
        \filldraw[red] (22,5.3)  circle (4pt);
        \filldraw[red] (23,6.3)  circle (4pt);
        \filldraw[red] (24,5.3)  circle (4pt);
        \filldraw[red] (25,4.3)  circle (4pt);
        \filldraw[red] (26,3.3)  circle (4pt);
        \filldraw[red] (27,2.3)  circle (4pt);
        \filldraw[red] (28,3.3)  circle (4pt);
        \filldraw[red] (29,2.3)  circle (4pt);
        \filldraw[red] (30,1.3)  circle (4pt);
        \filldraw[red] (31,0.3)  circle (4pt);
        \draw (13.1,12) node {$R$};
        \draw (15.5,1) node {Structure ($3b$)};
        \draw[solid,line width=0.2mm,<->] (23.5,6.5)--++(3.7,-3.7);
        \draw[solid,line width=0.2mm,<->] (17.5,11.5)--++(3.7,-3.7);
        \draw (20,10) node {$p$};
        \draw (26,5) node {$p$};
    \end{tikzpicture}\\
    \centering\begin{tikzpicture}[scale=0.4]
        \draw[solid,line width=0.4mm,blue] (0,0)-- ++\u -- ++\u -- ++\d -- ++\d -- ++\u -- ++\d -- ++\u -- ++\d;
        \filldraw[blue] (0,0)  circle (4pt);
        \filldraw[blue] (1,1)  circle (4pt);
        \filldraw[blue] (2,2)  circle (4pt);
        \filldraw[blue] (3,1)  circle (4pt);
        \filldraw[blue] (4,0)  circle (4pt);
        \filldraw[blue] (5,1)  circle (4pt);
        \filldraw[blue] (6,0)  circle (4pt);
        \filldraw[blue] (7,1)  circle (4pt);
        \filldraw[blue] (8,0)  circle (4pt);
         \draw[dashed] (0,0)--(8,0);
        \draw[solid,line width=0.4mm,red] (0,0.3)-- ++\u -- ++\u -- ++\u -- ++\d -- ++\d -- ++\u -- ++\d -- ++\d;
        \filldraw[red] (0,0.3)  circle (4pt);
        \filldraw[red] (1,1.3)  circle (4pt);
        \filldraw[red] (2,2.3)  circle (4pt);
        \filldraw[red] (3,3.3)  circle (4pt);
        \filldraw[red] (4,2.3)  circle (4pt);
        \filldraw[red] (5,1.3)  circle (4pt);
        \filldraw[red] (6,2.3)  circle (4pt);
        \filldraw[red] (7,1.3)  circle (4pt);
        \filldraw[red] (8,0.3)  circle (4pt);
        \draw (4.3,-1.5) node {Structure ($3c$)};
    \end{tikzpicture}\qquad
    \begin{tikzpicture}[scale=0.4]
        \draw[solid,line width=0.4mm,blue] (0,0)-- ++\u -- ++\d -- ++\u -- ++\d -- ++\u -- ++\d -- ++\u -- ++\d;
        \filldraw[blue] (0,0)  circle (4pt);
        \filldraw[blue] (1,1)  circle (4pt);
        \filldraw[blue] (2,0)  circle (4pt);
        \filldraw[blue] (3,1)  circle (4pt);
        \filldraw[blue] (4,0)  circle (4pt);
        \filldraw[blue] (5,1)  circle (4pt);
        \filldraw[blue] (6,0)  circle (4pt);
        \filldraw[blue] (7,1)  circle (4pt);
        \filldraw[blue] (8,0)  circle (4pt);
        \draw[dashed] (0,0)--(8,0);
        \draw[solid,line width=0.4mm,red] (0,0.3)-- ++\u -- ++\u -- ++\u -- ++\d -- ++\d -- ++\u -- ++\d -- ++\d;
        \filldraw[red] (0,0.3)  circle (4pt);
        \filldraw[red] (1,1.3)  circle (4pt);
        \filldraw[red] (2,2.3)  circle (4pt);
        \filldraw[red] (3,3.3)  circle (4pt);
        \filldraw[red] (4,2.3)  circle (4pt);
        \filldraw[red] (5,1.3)  circle (4pt);
        \filldraw[red] (6,2.3)  circle (4pt);
        \filldraw[red] (7,1.3)  circle (4pt);
        \filldraw[red] (8,0.3)  circle (4pt);
        \draw (4.3,-1.5) node {Structure ($3c$)};
    \end{tikzpicture}
    \caption{The structure of linear intervals $[P,Q]$ in $\FF_n^p$ when  $\texttt{type}(Q)=\texttt{type}(P)+1$. See Lemma~\ref{structure lin int F_n^p j-i=1}}
    \label{fig:strucure lin int j-i=1}
\end{figure}

\begin{proof}
    We can easily check that the three statements hold for $n\leq 4$ (see Figure~\ref{fig:strucure lin int j-i=1}). Thus, we suppose $n\geq 5$ which rules out the case ($c$). An interval of the form ($a$) or ($b$) is clearly linear, thus we focus on the converse.
    
    Let $[P,Q]$ be a linear interval. Since $j-i=1$, $P$ ends with $DUD^i$, and $Q$ ends with $UD^{i+1}$. \begin{enumerate}
        \item[(i)] Assume that $P$ ends with $UDUD^i$. If $[P,Q]=[U^iDUD^i,U^{i+1}D^{i+1}]$, then $[P,Q]$ has the form ($a$) with $k=1$ and $R$ empty. Otherwise, let $k$ (resp. $\ell$) be the greatest integer such that $P=P'(DU)^kD^i$ (resp. $\ell\leq k$ and $Q=Q'(UD)^\ell D^i$). For a contradiction, let us assume $\ell<k$. If $P'=U^i$ (this implies $k-\ell\geq2$), the paths $U^i(UD)^{k-\ell+1}(DU)^{\ell-1}D^i$ and $U^{i+2}D^2(UD)^{k-\ell-2}(DU)^\ell D^i$ both belong to $[P,Q]$, and they are not comparable, which contradicts the linearity of $[P,Q]$.
        If $P'=P''D$, then we obtain a contradiction with $P''D(UD)^{k-\ell+1}(DU)^{\ell-1}D^i$ and $P''UD^2(UD)^{k-\ell-1}(DU)^\ell D^i$. This proves that we necessarily have $k=\ell$ and thus $P=U^iP'(DU)^kD^i$ and $Q=U^iQ'(UD)^kD^i$ for some $P',Q'\in\F_{n-k-i}^p$ with types between 1 and $p-1$. We conclude that $P'=Q'$ as at the end of the proof of Lemma \ref{structure lin int F_n^p j-i > 2}. Finally we have that $[P,Q]$ has the form ($a$).
        \item[(ii)] Assume that $P$ ends with $D^2UD^i$. We have $P=P'DUD^i$ and $Q=Q'UD^{i+1}$. If $P'=Q'$, then $[P,Q]$ has the form ($a$) with $k=1$. So we assume $P'\ne Q'$. Let $m$, $1\leq m\leq p$ be the integer such that $P$ ends with $UD^mUD^i$, and $m-1\leq \ell\leq p$ such that $Q$ ends with $UD^\ell UD^{i+1}$. For a contradiction suppose $\ell=m-1$. Then $P=P''DUD^mUD^i$ and $Q=Q''DUD^{m-1}UD^{i+1}$. Since $P'\ne Q'$, there exists some $T\leq Q''$ and such that $T$ covers $P''$. So $TDUD^mUD^i$ and $P''DUD^{m-1}UD^{i+1}$ are not comparable which contradicts the linearity of $[P,Q]$. So, let us consider $\ell\geq m$. We cannot have $m<p$, since $P''DUD^{m-1}UD^{i+1}$ and $P''UD^{m+1}UD^i$ are not comparable which would contradict the linearity of $[P,Q]$. We conclude that $m=\ell=p$. Let $k$ be the greatest integer such that $P=P''D(UD^p)^kUD^i$ and $Q=Q''(UD^p)^kUD^{i+1}$. Using the maximality of $k$ and a similar argument as before, we can easily conclude that $P''=Q''$, proving that $[P,Q]$ has the form ($b$).
    \end{enumerate}
\end{proof}


Let $L_p(x,y)$ be the generating function where the coefficient of $x^ny^k$ in its series expansion is the number of  linear intervals of height $k$ in $\FF_n^p$. We handle separately the three cases $p=2$, $p\geq 3$, and $p=\infty$. 


\subsection{In the lattices $\FF_n^2$}
In this part, we fix $p=2$.
\begin{theorem}
    The generating function $L_2(x,y)$ of the number of linear intervals in $\FF_n^2$ with respect to $n$ and the interval height is given by 
    $$L_2(x,y)=\frac{x^{4} y^{4}+y^{3} x^{4}+1}{1-x -x^2}+\frac{x^{2} y \left(x^{2}-1\right) \left(x^{3} y^{2}-1\right)}{\left(x y -1\right) \left(x^{2}+x -1\right)^{2} \left(x^{2} y -1\right)}.$$
\end{theorem}
\begin{proof}
    Let $[P,Q]$ be a linear interval in $\FF_n^2$. According to the lemmas above, we have two cases to consider:
   \textbf{(i)} $\texttt{type}(P)=\texttt{type}(Q)$ and \textbf{(ii)} $\texttt{type}(P)=1$ and $\texttt{type}(Q)=~2$. Since $p=2$, Lemma~\ref{structure lin int F_n^p j-i > 2} is not considered.
   
    \textbf{Case (i)}. Using Lemma~\ref{structure lin int F_n^p j=i}, we have $[P,Q]=[P'UD,Q'UD]$ (resp. $[P,Q]=[UP'UD^2,UQ'UD^2]$) where $[P',Q']$ is a linear interval in $\F_{n-1}^2$ (resp. $\F_{n-2}^2$). Then the contribution of these intervals is $$xL_2(x,y)+x^2L_2(x,y).$$ 
    
    \textbf{Case (ii)}. Using Lemma~\ref{structure lin int F_n^p j-i=1}, the contribution of intervals from the statement (3$a$) is $$\frac{x^2y}{1-xy}\cdot (xF(x)+1),$$ where $xF(x)+1$ is the generating function for the empty path and paths in $\F^2$ of type 1, with $F(x)=\frac{1}{1-x-x^2}$ is the generating function for all paths in $\F^2$. The contribution of intervals with structure (3$b$) is $$\frac{x^4y^2}{1-x^2y}\cdot (xF(x)+1).$$

    The contribution of intervals with structure (3$c$) is $x^4y^3+x^4y^4$. Finally, we obtain the following equation for $L_2(x,y)$:
   $$L_2(x,y)=1+x(1+x)L_2(x,y)+x^4y^3+x^4y^4+\left(\frac{x^2y}{1-xy}+ \frac{x^4y^2}{1-x^2y}   \right) \cdot (xF(x)+1), $$
    which induces the desired expression of $L_2(x,y)$.
\end{proof}

The first terms of the series expansion of $L_2(x,y)$ are
\begin{align*}1 + x + (2+y)x^2 + (3+2y+y^2)x^3 + (5+4y+3y^2+2y^3+y^4)x^4 \\+ (8+8y+6y^2+3y^3+2y^4)x^5 + (13+15y+12y^2+7y^3+4y^4+y^5)x^6+O(x^7).
\end{align*}

\begin{corollary}
The generating function of the number of linear intervals in $\FF_n^2$ is given by 
    $$L_2(x,1)=\frac{1-x+x^3+3x^4-2x^5-2x^6}{\left(x^{2}+x -1\right)^{2}},$$
    and the coefficient of $x^n$ in the series expansion is $[x^0]L_2(x,1)=[x^1]L_2(x,1)=1$, $[x^2]L_2(x,1)=3$, and
     $$[x^n]L_2(x,1)=\frac{4 \left(n +5\right) \mathit{F}_n-\left(2 n +27\right) \mathit{F}_{n-1}}{5} \ \text{ for } \ n\geq3,$$
    where $F_n$ is the $n$-th Fibonacci number (see Section \ref{sec:intro}). An asymptotic approximation for the number of linear intervals in $\FF_n^2$ is $$\frac{3n}{5}\left(\frac{1+\sqrt{5}}{2}\right)^n.$$
\end{corollary}
\begin{proof} The closed form of the $n$-th term is obtain by decomposing the rational fraction into partial fractions. The asymptotic is easily obtained  from a classical singularity analysis, see for instance \cite{analytic combinatorics}.
\end{proof}

\begin{corollary}
    The generating function of the limit distribution of the height of the linear intervals in $\FF_n^2$ as $n\to\infty$ is given by
    $$\frac{y((7-3\sqrt{5})y^2 + (7-3\sqrt{5})y -2)}{6((2 - \sqrt{5})y^2 + y - 1)}.$$
    In particular it has expected value $\frac{3+\sqrt{5}}{2}$, and variance $\frac{7+2\sqrt{5}}{3}$.
\end{corollary}

\begin{proof}
     We use singularity analysis, see \cite{analytic combinatorics}. The main singularity in $x$ of $L_2(x,y)$ is $\frac{\sqrt{5}-1}{2}$, in particular it does not depend on $y$. The limit law of the height of the linear intervals in $\F_n^2$ as $n\to\infty$ is then discrete. Near this singularity, we have the approximation
    $$L_2(x,y)\sim \frac{(-16y^3 - 16y^2 + 2y)\sqrt{5} + 36y(y^2 + y - 1/6)}{5(y\sqrt{5} - y - 2)(y\sqrt{5} - 3y + 2)}\left(x-\frac{\sqrt{5}-1}{2}\right)^{-2}.$$
    We deduce that
    $$\frac{[x^n]L_2(x,y)}{[x^n]L_2(x,1)}\underset{n\to\infty}{\longrightarrow}\frac{(-16y^3 - 16y^2 + 2y)\sqrt{5} + 36y(y^2 + y - 1/6)}{5(y\sqrt{5} - y - 2)(y\sqrt{5} - 3y + 2)}\cdot\frac{5}{2}\cdot\frac{3+\sqrt{5}}{2},$$
    which yields the desired generating function after simplification.
\end{proof}

\subsection{In the lattices $\FF_n^p$, $p\geq3$} In this part, we fix $p\geq 3$. Recall that the generating function $F_p(x)$ for the number of elements in $\FF_n^p$, $n\geq 0$, is $F_p(x)=\frac{1}{1-x-x^2-\cdots -x^p}$. We also set $$G_p(x):=1-x-x^2-\ldots-x^p.$$

\begin{thm} The generating function $L_p(x,y)$ of the number of linear intervals in $\FF_n^p$, $n\geq 0$, with respect to $n$ and the interval height is 
     $$L_p(x,y)=F_p(x)\cdot\left(1+\frac{x^3y^3(1-x^{p-2})}{1-x}+V_p(x,y)+W_p(x,y)\right),$$
    where 
    \begin{align*}
        V_p(x,y)&=\sum_{\substack{1\leq i,j\leq p \\ j-i\geq2}}x^jy^{j-i}\left(1+\frac{x+x^2+\ldots+x^{p-(j-i)}}{G_p(x)}\right),\\
    W_p(x,y)&=
    \frac{y(1-x^p)(x^2-x^{p+1})(1-x^{p+1}y^2)}{(1-x)(1-xy)(1-x^py)G_p(x)}.\end{align*}
\end{thm}
\begin{proof}
   Let $[P,Q]$ be a linear interval in $\FF_n^p$. According to the lemmas above, we have three cases to consider:
   (i) $\texttt{type}(P)=\texttt{type}(Q)$, (ii) $\texttt{type}(Q)\geq\texttt{type}(Q)+2$, and (iii) $\texttt{type}(Q)=\texttt{type}(Q)+1$.
   
    \textbf{Case (i)}. Using Lemma~\ref{structure lin int F_n^p j=i}, we have $[P,Q]=[U^{i-1}P'UD^i,U^{i-1}Q'UD^i]$, $1\leq i\leq p$, where $[P',Q']$ is a linear interval in $\FF_{n-i}^p$. Then the contribution of these intervals is $$xL_p(x,y)+x^2L_p(x,y)+\ldots+x^pL_p(x,y)=(1-G_p(x))\cdot L_p(x,y).$$

    \textbf{Case (ii)}. Using Lemma~\ref{structure lin int F_n^p j-i > 2}, the contribution of paths of structure (2$a$) is given by 
    $$x^3y^3\sum_{i=0}^{p-3}x^i=\frac{x^3y^3(1-x^{p-2})}{1-x}.$$

   Let $V_p(x,y)$ be the contribution of paths of  structure (2$b$). Since $R$ is either empty or has type between $1$ and $p-(j-i)$, the generating function for these paths $R$ is 
   $$1+(x+x^2+\ldots+x^{p-(j-i)})F_p(x).$$ 
   Multiplying by $x^jy^{j-i}$ and summing for  $1\leq i,j\leq p$ with $j-i\geq 2$, we obtain the desired result for $V_p(x,y)$:
  
         $$ V_p(x,y)=\sum_{\substack{1\leq i,j\leq p \\ j-i\geq2}}x^jy^{j-i}\left(1+\frac{x+x^2+\ldots+x^{p-(j-i)}}{G_p(x)}\right).$$

   \textbf{Case (iii)}. Using Lemma~\ref{structure lin int F_n^p j-i=1}, the contribution for paths belonging of structure (3$a$) is the contribution of $R$ of type in $[1,p-1]$ (\textit{i.e.} $1+(x+x^2+\ldots+x^{p-1})F_p(x)$) multiplied by the contribution of $[U^i(DU)^kD^i, U^i(UD)^kD^i]$ for $k\geq 1$ and $1\leq i\leq p-1$, which gives 
   $$\sum_{k=1}^{+\infty} x^{i+k}y^k=\frac{x^{i+1}y}{1-xy}.$$

    For structure (3$b$), $R$ has the same contribution as for the previous case (3$a$), and we need to multiply by the contribution of     $[U^{k(p-1)+i}D(UD^p)^kUD^i,U^{k(p-1)+i}(UD^p)^kUD^{i+1}]$, with $k\geq1$ and $1\leq i\leq p-1$.
    So the generating function is 
    $$\sum_{i=1}^{p-1}\sum_{k=1}^{+\infty} x^{i+1+kp}y^{k+1}=\frac{x^{p+i+1}y^2}{1-x^py}.$$
    Paths of structures (3$c$) have no contribution for $p\geq 3$.
    
    Finally, the generating function for structure $(3)$ is 
    $$W_p(x,y)=\sum_{i=1}^{p-1}\left(\frac{x^{i+1}y}{1-xy}+\frac{x^{p+i+1}y^2}{1-x^py}\right)\left(1+\frac{x+x^2+\ldots+x^{p-1}}{G_p(x)}\right),$$
    which can be simplified as in the Theorem after computation. Finally, adding everything up we obtain the expected result for $L_p(x,y)$. 
\end{proof}

\begin{remark} We did not succeed to obtain a nice closed form for $L_p(x,y)$. However,  even though the expression of $L_p(x,y)$ is quite heavy, we can observe that its main singularity is the smallest root $r_p$ of $G_p(x)$. In particular it does not depend on $y$. This yields, see for instance \cite{analytic combinatorics}, that the limit law of the length of the linear intervals in $\FF_n^p$ as $n\to\infty$ is discrete.  Furthermore, $r_p$ is a singularity of multiplicity 2 in $L_p(x,1)$, thus the number of linear intervals in $\FF_n^p$ is asymptotically  $c\cdot n\cdot r_p^{-n}$ for some constant $c$ (see \cite{analytic combinatorics}), so by Corollary \ref{nb coverings}, it is proportional to the number of coverings.
\end{remark}


\subsection{In the lattice $\FF_n^\infty$}
The generating function $L_\infty(x,y)$ for the linear interval in  $\FF_n^\infty$ is obtained using a discrete continuity argument.

\begin{corollary} The generating function of the linear intervals in $\FF_n^\infty$ with respect to $n$ and the interval height is given by
    $$L_\infty(x,y)=\frac{1-y^2(1+y)^2x^4+2x^5y^4-(3-y-y^2)x^3y+2(2y+1)x^2-(3+y)x}{(1-xy)(1-2x)^2}.$$
\end{corollary}
\begin{proof}
    Once again we could go through the structure of the linear intervals in $\FF_n^\infty$ and do a similar proof as for $\FF_n^p$, but we can also use a discrete continuity argument, as in Corollary \ref{boolean F infty}. Indeed, with the same argument, we have $L_p\to L_\infty$ as $p\to\infty$. We have the following limits as $p\to\infty$: $$G_p(x)\to\frac{1-2x}{1-x},\quad\frac{x^3y^3(1-x^{p-2})}{1-x}\to \frac{x^3y^3}{1-x}, \quad W_p(x,y)\to\frac{x^2y}{(1-xy)(1-2x)},$$ and with the help of a computer algebra program, $V_p(x,y)\to\frac{x^3y^2}{(1-xy)(1-2x)}$. We deduce that 
    $$L_\infty(x,y)=\frac{1-x}{1-2x}\left(1+\frac{x^3y^3}{1-x}+\frac{x^2y}{(1-xy)(1-2x)}+\frac{x^3y^2}{(1-xy)(1-2x)}\right),$$
    which gives the desired expression.
\end{proof}

\begin{corollary} The generating function for the number of linear intervals in $\FF_n^\infty$
is $$L_\infty(x,1)={\frac{1-3x+3x^2+2x^3-2x^4}{ \left( 1-2x
 \right) ^{2}}},$$
  and the coefficient of $x^n$ in the series expansion is $[x^0]L_\infty(x,1)=[x^1]L_\infty(x,1)=1$, $[x^2]L_\infty(x,1)=3$, and
  $$[x^n]L_\infty(x,1)=(3n+1)\cdot2^{n-3} \ \text{ for } n\geq3.$$
\end{corollary}

\begin{corollary}
    The generating function of the limit distribution of the length of the linear intervals in $\FF_n^\infty$ as $n\to\infty$ is given by
    $$\frac{y(2+y)}{3(2-y)}.$$
    This is a `geometric-like' law with parameter $1/2$. It has expected value $7/3$, variance $20/9$, and if $p_k$ denotes the asymptotic proportion of linear intervals having length $k$, then
    $$\left\{\begin{array}{ccl}
        p_0 &=& 0,  \\
        p_1 &=& 1/3, \\
        p_k &=& \frac{1}{3\cdot2^{k-2}} \text{ for } k\geq2.
    \end{array}\right.$$
\end{corollary}
\begin{proof}
    Near $1/2$, the main singularity of $L_\infty(x,y)$, we have the following approximation:
    $$L_\infty(x,y)\sim\frac{y(2+y)}{32(2-y)}\left(\frac{1}{2}-x\right)^{-2}.$$
    We then easily deduce the desired generating function:
    $$\frac{[x^n]L(x,y)}{[x^n]L(x,1)}\underset{n\to+\infty}{\longrightarrow}\frac{y(2+y)}{3(2-y)}=\frac{y}{3}+\sum_{k\geq 2}\frac{y^k}{3\cdot2^{k-2}}.$$
\end{proof}

\section{Intervals}\label{sec:intervals}

In this section we count intervals in $\FF_n^2$ and $\FF_n^\infty$. 
A crucial point in our study is the following three facts which can be checked with a simple observation.

\begin{fact}\label{fac:first-insert}
  Let $m, n\ge 1$, and let $[P',Q']$ be an interval of $\FF_{n}^2$ (resp. $\FF_n^\infty$). Delete in $P'$ and $Q'$ the first peak $UD$, to obtain two paths $P$ and $Q$ of $\FF_{n-1}^2$ (resp. $\FF_{n-1}^\infty$). Then they form an interval $[P,Q]$.
\end{fact}

\begin{fact}\label{fac:first-insertconv1}  Conversely, start from an interval $[P,Q]=[U^iD^k\alpha,U^jD^{\ell}\beta]$ in $\FF_{n-1}^\infty$ where $k,\ell\geq1$ are maximal and $\alpha$, $\beta$ possibly empty.  Inserting a  peak $UD$ in the first ascents of $P$ and $Q$, starting at heights $a\in [ 0, i] $ and $b\in [0, j]$ respectively, yields an interval $[P',Q']$ in $\FF_n^\infty$ if and only if $a\in\{i-1,i\}$, $b\in\{j-1,j\}$ and $a\leq b$.
\end{fact}

\begin{fact}\label{fac:first-insertconv2}  Conversely, start from an interval $[P,Q]=[U^iD^k\alpha,U^jD^{\ell}\beta]$ in $\FF_{n-1}^2$ where $k,\ell\in\{1,2\}$ are maximal and $\alpha$, $\beta$ possibly empty.  Inserting a  peak $UD$ in the first ascents of $P$ and $Q$, starting at heights $a\in  [ 0, i] $ and $b\in [0, j]$ respectively, yields an interval $[P',Q']$ in $\FF_{n}^2$ if and only if 
\begin{itemize}
    \item[(1)]  $a\in\{i-1,i\}$, $b\in\{j-1,j\}$, $a\leq b$, when $k=\ell=1$,
    \item[(2)]  $a\in\{i-1,i\}$, $b=j-1$, $a\leq b$, when $k=1$, $\ell=2$,
    \item[(3)]  $a=i-1$, $b\in\{j-1,j\}$, when $k=2$, $\ell=1$,
     \item[(4)]  $a=i-1$, $b=j-1$, when $k=2$, $\ell=2$.
    \end{itemize}
\end{fact}

\subsection{In the lattice $\FF_n^\infty$} According to Fact~\ref{fac:first-insertconv1}, the rule that describes the first ascent lengths of paths $P'\le Q'$ obtained from  $[P,Q]\in\FF_{n}^\infty$ in terms of the first ascent lengths $a$ and $b$ of $P$ and $Q$ is
\begin{equation}
  (a,b)  \rightarrow \left\{
    \begin{array}{ll}
    (a,b+1), (a+1,b+1), (a,b), & \text{ if } b-a=0,  \\
    (a,b+1), (a+1,b+1), (a,b), (a+1,b), & \text{ if } b-a\geq 1,  \\
       \end{array}
    \right.
    \label{gen rules motzkin}
\end{equation}
starting with the root $(1,1)$ corresponding to the interval $[UD,UD]$. In fact, the above rules describe the construction of some  paths in the first quadrant of the plane.

\begin{theorem}\label{cor:bij-mn} There is a bijection between  intervals in $\FF^\infty_{n}$ and bicolored Motzkin paths of length $n-1$ in the quarter plane, \textit{i.e.},  paths in the quarter plane starting at $(0,0)$ and consisting  of $n-1$ steps $\U=(1,1)$, $\D=(1,-1)$, $\Fb=(1,0)$, $\Fg=(1,0)$.
\end{theorem}

\begin{figure}[h]\label{thrule1}
    \centering
    \captionsetup{justification=centering,margin=0.5cm}
    \begin{tikzpicture}
    \node (1) at (0,0) {\begin{tikzpicture}[scale=0.3]
        \draw[solid,line width=0.4mm,blue] (0,0)-- ++\u -- ++\d;
        \draw[solid,line width=0.4mm,red] (0,0.6)-- ++\u -- ++\d;
        \filldraw[blue] (0,0)  circle (4pt);
        \filldraw[blue] (1,1)  circle (4pt);
        \filldraw[blue] (2,0)  circle (4pt);
        \filldraw[red] (0,0.5)  circle (4pt);
        \filldraw[red] (1,1.5)  circle (4pt);
        \filldraw[red] (2,0.5)  circle (4pt);
    \end{tikzpicture}};
    \node (2) at (2,0) {\begin{tikzpicture}[scale=0.3]
        \draw[solid,line width=0.4mm,blue] (0,0)-- ++\u -- ++\d-- ++\u -- ++\d;
        \draw[solid,line width=0.4mm,red] (0,0.6)-- ++\u-- ++\u -- ++\d -- ++\d;
        \filldraw[blue] (0,0)  circle (4pt);
        \filldraw[blue] (1,1)  circle (4pt);
        \filldraw[blue] (2,0)  circle (4pt);
        \filldraw[blue] (3,1)  circle (4pt);
        \filldraw[blue] (4,0)  circle (4pt);
        \filldraw[red] (0,0.5)  circle (4pt);
        \filldraw[red] (1,1.5)  circle (4pt);
        \filldraw[red] (2,2.5)  circle (4pt);
        \filldraw[red] (3,1.5)  circle (4pt);
        \filldraw[red] (4,0.5)  circle (4pt);
    \end{tikzpicture}};
    \node (3) at (4.6,0) {\begin{tikzpicture}[scale=0.3]
        \draw[solid,line width=0.4mm,blue] (0,0)-- ++\u -- ++\d-- ++\u -- ++\d-- ++\u -- ++\d;
        \draw[solid,line width=0.4mm,red] (0,0.6)-- ++\u-- ++\u -- ++\d-- ++\u  -- ++\d -- ++\d;
        \filldraw[blue] (0,0)  circle (4pt);
        \filldraw[blue] (1,1)  circle (4pt);
        \filldraw[blue] (2,0)  circle (4pt);
        \filldraw[blue] (3,1)  circle (4pt);
        \filldraw[blue] (4,0)  circle (4pt);
        \filldraw[blue] (5,1)  circle (4pt);
        \filldraw[blue] (6,0)  circle (4pt);
        \filldraw[red] (0,0.5)  circle (4pt);
        \filldraw[red] (1,1.5)  circle (4pt);
        \filldraw[red] (2,2.5)  circle (4pt);
        \filldraw[red] (3,1.5)  circle (4pt);
        \filldraw[red] (4,2.5)  circle (4pt);
        \filldraw[red] (5,1.5)  circle (4pt);
        \filldraw[red] (6,0.5)  circle (4pt);
    \end{tikzpicture}};
    \node (4) at (8,0) {\begin{tikzpicture}[scale=0.3]
        \draw[solid,line width=0.4mm,blue] (0,0)-- ++\u -- ++\d-- ++\u -- ++\d-- ++\u -- ++\d-- ++\u -- ++\d;
        \draw[solid,line width=0.4mm,red] (0,0.6)-- ++\u-- ++\u-- ++\u -- ++\d -- ++\d-- ++\u  -- ++\d -- ++\d;
        \filldraw[blue] (0,0)  circle (4pt);
        \filldraw[blue] (1,1)  circle (4pt);
        \filldraw[blue] (2,0)  circle (4pt);
        \filldraw[blue] (3,1)  circle (4pt);
        \filldraw[blue] (4,0)  circle (4pt);
        \filldraw[blue] (5,1)  circle (4pt);
        \filldraw[blue] (6,0)  circle (4pt);
        \filldraw[blue] (7,1)  circle (4pt);
        \filldraw[blue] (8,0)  circle (4pt);
        \filldraw[red] (0,0.5)  circle (4pt);
        \filldraw[red] (1,1.5)  circle (4pt);
        \filldraw[red] (2,2.5)  circle (4pt);
        \filldraw[red] (3,3.5)  circle (4pt);
        \filldraw[red] (4,2.5)  circle (4pt);
        \filldraw[red] (5,1.5)  circle (4pt);
        \filldraw[red] (6,2.5)  circle (4pt);
        \filldraw[red] (7,1.5)  circle (4pt);
        \filldraw[red] (8,0.5)  circle (4pt);
    \end{tikzpicture}};
    \node (5) at (12,0) {\begin{tikzpicture}[scale=0.3]
        \draw[solid,line width=0.4mm,blue] (0,0)-- ++\u-- ++\u -- ++\d -- ++\d-- ++\u -- ++\d-- ++\u -- ++\d-- ++\u -- ++\d;
        \draw[solid,line width=0.4mm,red] (0,0.6)-- ++\u-- ++\u-- ++\u -- ++\d-- ++\u -- ++\d -- ++\d-- ++\u  -- ++\d -- ++\d;
        \filldraw[blue] (0,0)  circle (4pt);
        \filldraw[blue] (1,1)  circle (4pt);
        \filldraw[blue] (2,2)  circle (4pt);
        \filldraw[blue] (3,1)  circle (4pt);
        \filldraw[blue] (4,0)  circle (4pt);
        \filldraw[blue] (5,1)  circle (4pt);
        \filldraw[blue] (6,0)  circle (4pt);
        \filldraw[blue] (7,1)  circle (4pt);
        \filldraw[blue] (8,0)  circle (4pt);
        \filldraw[blue] (9,1)  circle (4pt);
        \filldraw[blue] (10,0)  circle (4pt);
        \filldraw[red] (0,0.5)  circle (4pt);
        \filldraw[red] (1,1.5)  circle (4pt);
        \filldraw[red] (2,2.5)  circle (4pt);
        \filldraw[red] (3,3.5)  circle (4pt);
        \filldraw[red] (4,2.5)  circle (4pt);
        \filldraw[red] (5,3.5)  circle (4pt);
        \filldraw[red] (6,2.5)  circle (4pt);
        \filldraw[red] (7,1.5)  circle (4pt);
        \filldraw[red] (8,2.5)  circle (4pt);
        \filldraw[red] (9,1.5)  circle (4pt);
        \filldraw[red] (10,0.5)  circle (4pt);
    \end{tikzpicture}};
    \node (6) at (12.5,-2) {\begin{tikzpicture}[scale=0.3]
        \draw[solid,line width=0.4mm,blue] (0,0)-- ++\u-- ++\u-- ++\u -- ++\d -- ++\d -- ++\d-- ++\u -- ++\d-- ++\u -- ++\d-- ++\u -- ++\d;
        \draw[solid,line width=0.4mm,red] (0,0.6)-- ++\u-- ++\u-- ++\u-- ++\u -- ++\d -- ++\d-- ++\u -- ++\d -- ++\d-- ++\u  -- ++\d -- ++\d;
        \filldraw[blue] (0,0)  circle (4pt);
        \filldraw[blue] (1,1)  circle (4pt);
        \filldraw[blue] (2,2)  circle (4pt);
        \filldraw[blue] (3,3)  circle (4pt);
        \filldraw[blue] (4,2)  circle (4pt);
        \filldraw[blue] (5,1)  circle (4pt);
        \filldraw[blue] (6,0)  circle (4pt);
        \filldraw[blue] (7,1)  circle (4pt);
        \filldraw[blue] (8,0)  circle (4pt);
        \filldraw[blue] (9,1)  circle (4pt);
        \filldraw[blue] (10,0)  circle (4pt);
        \filldraw[blue] (11,1)  circle (4pt);
        \filldraw[blue] (12,0)  circle (4pt);
        \filldraw[red] (0,0.5)  circle (4pt);
        \filldraw[red] (1,1.5)  circle (4pt);
        \filldraw[red] (2,2.5)  circle (4pt);
        \filldraw[red] (3,3.5)  circle (4pt);
        \filldraw[red] (4,4.5)  circle (4pt);
        \filldraw[red] (5,3.5)  circle (4pt);
        \filldraw[red] (6,2.5)  circle (4pt);
        \filldraw[red] (7,3.5)  circle (4pt);
        \filldraw[red] (8,2.5)  circle (4pt);
        \filldraw[red] (9,1.5)  circle (4pt);
        \filldraw[red] (10,2.5)  circle (4pt);
        \filldraw[red] (11,1.5)  circle (4pt);
        \filldraw[red] (12,0.5)  circle (4pt);
    \end{tikzpicture}};
    \node (7) at (7,-2) {\begin{tikzpicture}[scale=0.3]
        \draw[solid,line width=0.4mm,blue] (0,0)-- ++\u-- ++\u-- ++\u -- ++\d-- ++\u -- ++\d -- ++\d -- ++\d-- ++\u -- ++\d-- ++\u -- ++\d-- ++\u -- ++\d;
        \draw[solid,line width=0.4mm,red] (0,0.6)-- ++\u-- ++\u-- ++\u-- ++\u-- ++\u -- ++\d -- ++\d -- ++\d-- ++\u -- ++\d -- ++\d-- ++\u  -- ++\d -- ++\d;
        \filldraw[blue] (0,0)  circle (4pt);
        \filldraw[blue] (1,1)  circle (4pt);
        \filldraw[blue] (2,2)  circle (4pt);
        \filldraw[blue] (3,3)  circle (4pt);
        \filldraw[blue] (4,2)  circle (4pt);
        \filldraw[blue] (5,3)  circle (4pt);
        \filldraw[blue] (6,2)  circle (4pt);
        \filldraw[blue] (7,1)  circle (4pt);
        \filldraw[blue] (8,0)  circle (4pt);
        \filldraw[blue] (9,1)  circle (4pt);
        \filldraw[blue] (10,0)  circle (4pt);
        \filldraw[blue] (11,1)  circle (4pt);
        \filldraw[blue] (12,0)  circle (4pt);
        \filldraw[blue] (13,1)  circle (4pt);
        \filldraw[blue] (14,0)  circle (4pt);
        \filldraw[red] (0,0.5)  circle (4pt);
        \filldraw[red] (1,1.5)  circle (4pt);
        \filldraw[red] (2,2.5)  circle (4pt);
        \filldraw[red] (3,3.5)  circle (4pt);
        \filldraw[red] (4,4.5)  circle (4pt);
        \filldraw[red] (5,5.5)  circle (4pt);
        \filldraw[red] (6,4.5)  circle (4pt);
        \filldraw[red] (7,3.5)  circle (4pt);
        \filldraw[red] (8,2.5)  circle (4pt);
        \filldraw[red] (9,3.5)  circle (4pt);
        \filldraw[red] (10,2.5)  circle (4pt);
        \filldraw[red] (11,1.5)  circle (4pt);
        \filldraw[red] (12,2.5)  circle (4pt);
        \filldraw[red] (13,1.5)  circle (4pt);
        \filldraw[red] (14,0.5)  circle (4pt);
    \end{tikzpicture}};
    \node (8) at (1,-2) {\begin{tikzpicture}[scale=0.3]
        \draw[solid,line width=0.4mm,blue] (0,0)-- ++\u-- ++\u-- ++\u -- ++\d-- ++\u -- ++\d-- ++\u -- ++\d -- ++\d -- ++\d-- ++\u -- ++\d-- ++\u -- ++\d-- ++\u -- ++\d;
        \draw[solid,line width=0.4mm,red] (0,0.6)-- ++\u-- ++\u-- ++\u-- ++\u-- ++\u -- ++\d-- ++\u -- ++\d -- ++\d -- ++\d-- ++\u -- ++\d -- ++\d-- ++\u  -- ++\d -- ++\d;
        \filldraw[blue] (0,0)  circle (4pt);
        \filldraw[blue] (1,1)  circle (4pt);
        \filldraw[blue] (2,2)  circle (4pt);
        \filldraw[blue] (3,3)  circle (4pt);
        \filldraw[blue] (4,2)  circle (4pt);
        \filldraw[blue] (5,3)  circle (4pt);
        \filldraw[blue] (6,2)  circle (4pt);
        \filldraw[blue] (7,3)  circle (4pt);
        \filldraw[blue] (8,2)  circle (4pt);
        \filldraw[blue] (9,1)  circle (4pt);
        \filldraw[blue] (10,0)  circle (4pt);
        \filldraw[blue] (11,1)  circle (4pt);
        \filldraw[blue] (12,0)  circle (4pt);
        \filldraw[blue] (13,1)  circle (4pt);
        \filldraw[blue] (14,0)  circle (4pt);
        \filldraw[blue] (15,1)  circle (4pt);
        \filldraw[blue] (16,0)  circle (4pt);
        \filldraw[red] (0,0.5)  circle (4pt);
        \filldraw[red] (1,1.5)  circle (4pt);
        \filldraw[red] (2,2.5)  circle (4pt);
        \filldraw[red] (3,3.5)  circle (4pt);
        \filldraw[red] (4,4.5)  circle (4pt);
         \filldraw[red] (5,5.5)  circle (4pt);
         \filldraw[red] (6,4.5)  circle (4pt);
         \filldraw[red] (7,5.5)  circle (4pt);
         \filldraw[red] (8,4.5)  circle (4pt);
         \filldraw[red] (9,3.5)  circle (4pt);
         \filldraw[red] (10,2.5)  circle (4pt);
         \filldraw[red] (11,3.5)  circle (4pt);
         \filldraw[red] (12,2.5)  circle (4pt);
         \filldraw[red] (13,1.5)  circle (4pt);
        \filldraw[red] (14,2.5)  circle (4pt);
         \filldraw[red] (15,1.5)  circle (4pt);
         \filldraw[red] (16,0.5)  circle (4pt);
    \end{tikzpicture}};
    \draw[->] (1) -- (2);
    \draw[->] (2) -- (3);
    \draw[->] (3) -- (4);
    \draw[->] (4) -- (5);
    \draw[->] (14,0) arc (90:-60:1);
    \draw[->] (6) -- (7);
    \draw[->] (7) -- (8);
    \node (a) at (0.9,0.5) {$\U$};
    \node (b) at (3.2,0.5) {$\Fg$};
    \node (c) at (6.2,0.5) {$\U$};
    \node (d) at (9.9,0.5) {$\D$};
    \node (e) at (14.5,-1) {$\Fb$};
    \node (f) at (9.95,-1.5) {$\U$};
    \node (g) at (4.2,-1.5) {$\Fg$};
    \end{tikzpicture}
    \begin{tikzpicture}
        \draw[solid,line width=0.4mm] (-0.5,0)--(7.5,0);
        \draw[solid,line width=0.4mm] (0,-0.5)--(0,2.5);
        \draw[dashed,line width=0.2mm] (0,1)--(7.3,1);
        \draw[dashed,line width=0.2mm] (0,2)--(7.3,2);
        \draw[dashed,line width=0.2mm] (1,0)--(1,2.3);
        \draw[dashed,line width=0.2mm] (2,0)--(2,2.3);
        \draw[dashed,line width=0.2mm] (3,0)--(3,2.3);
        \draw[dashed,line width=0.2mm] (4,0)--(4,2.3);
        \draw[dashed,line width=0.2mm] (5,0)--(5,2.3);
        \draw[dashed,line width=0.2mm] (6,0)--(6,2.3);
        \draw[dashed,line width=0.2mm] (7,0)--(7,2.3);
        \draw[solid,line width=0.4mm,blue,->] (0,0)--(1,1);
        \draw[solid,line width=0.4mm,Green,->] (1,1)--(2,1);
        \draw[solid,line width=0.4mm,blue,->] (2,1)--(3,2);
        \draw[solid,line width=0.4mm,blue,->] (3,2)--(4,1);
        \draw[solid,line width=0.4mm,blue,->] (4,1)--(5,1);
        \draw[solid,line width=0.4mm,blue,->] (5,1)--(6,2);
        \draw[solid,line width=0.4mm,Green,->] (6,2)--(7,2);
    \end{tikzpicture}
\caption{The generation of the interval $[U^2(UD)^3D^2(UD)^3,U^4(UD)^2D^2(UD^2)^2]$ using the rules in the proof of Theorem~\ref{cor:bij-mn}. This interval is thus associated with the bicolored Motzkin path $\U\Fg\U\D\Fb\U\Fg$.}
\label{ex bij int binary motzkin}
\end{figure}
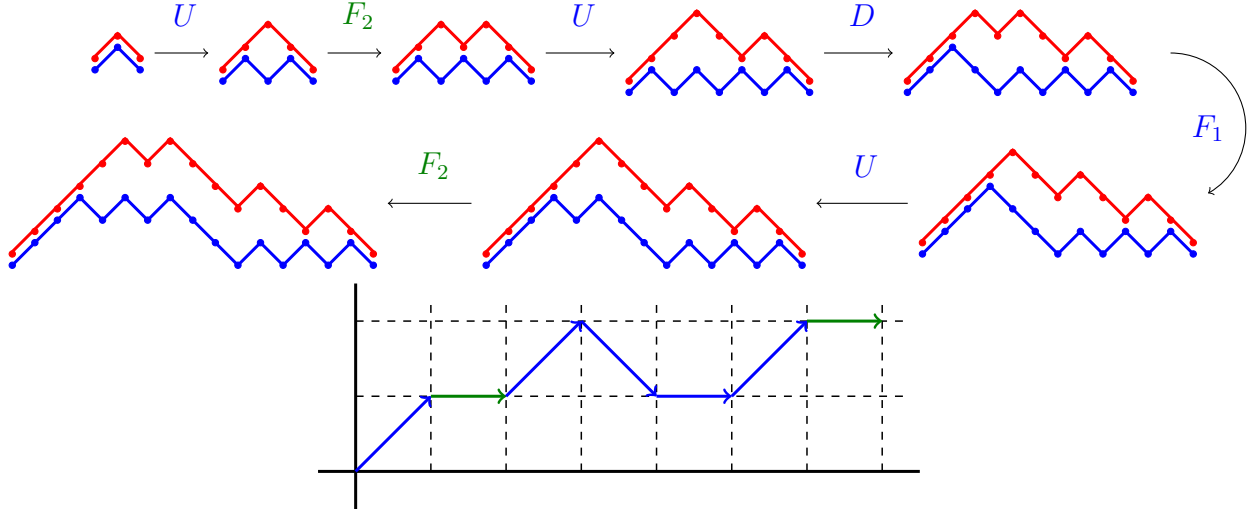

\begin{proof} To see the bijection, it suffices to rewrite the above rules in an equivalent system (in the sense where the rules generate an isomorphic generating tree). We start at the origin $(0,0)$ of the plane (corresponding to the interval $[UD,UD]$), and after $x$ steps, if we are at the point $(x,y)$, we can jump to the following points
$$
  (x,y)  \rightarrow \left\{
    \begin{array}{ll}
    (x+1,y+1), (x+1,y)_1, (x+1,y)_2, & \text{ if } y=0,  \\
    (x+1,y+1), (y+1,y)_1, (x+1,y)_2, (x+1,y-1), & \text{ if } y>0.  \\
       \end{array}
    \right.
    \label{rr}
  $$
   In comparison to (\ref{gen rules motzkin}), $x$ is the number of steps, and $y=b-a$. Note that we use the subscripts $1$ and $2$ to distinguish the step $(a,b)\to (a,b)$ from $(a,b)\to (a+1,b+1)$. The bicolored Motzkin path is then constructed from the origin $(0,0)$ using steps $\U=(1,1)$, $\D=(1,-1)$, $\Fb=(1,0)$, $\Fg=(1,0)$ corresponding to the rules $(x,y)\rightarrow (x+1,y+1)$, $(x,y)\rightarrow (x+1,y-1)$, $(x,y)\rightarrow (x+1,y)_1$ and $(x,y)\rightarrow (x+1,y)_2$, respectively. See Figure~\ref{ex bij int binary motzkin} for an illustration of the bijection. 
\end{proof}

\begin{thm}\label{int binary}
    The generating function $I(x,y)$  where the coefficient $x^ny^k$ is the number of intervals $[P,Q]\in\FF_n^\infty$ such that the difference between the lengths of the first ascent of $Q$ and $P$ equals $k$, is given by 
    $$I(x,y)=1+\frac{2x}{1-2x-2xy+\sqrt{1-4x}}.$$
    The generating function for the number of intervals is 
    $$I(x,1)
    =\frac{1}{2}\left(1+\frac{1}{\sqrt{1-4x}}\right),$$
    and the coefficient of $x^n$ in the series expansion is given by $\binom{2n-1}{n}$. (\seqnum{A001700} with a shift in \cite{oeis}).
\end{thm}
\begin{proof}
    By Theorem \ref{cor:bij-mn}, it suffices to count the bicolored Motzkin paths of length $n$ ending at ordinate $k$. Let $\mathcal{M}_n^k$ be the set of these paths and $\mathcal{M}^k$ the set of these paths of any length. Let $M(x,y)$ be the generating function for these paths, $x$ tracking the length and $y$ the final height.   A nonempty Motzkin path $M$ ending at height $k\geq 0$ can be decomposed either $M=\U M_1$ with $M_1\in \mathcal{M}_{n-1}^{k-1}$, or $M=\Fb M_1$ or $M=\Fg M_1$ with $M_1\in \mathcal{M}_{n-1}^{k}$, or  $\U M_1 \D M_2$ where $M_1\in \mathcal{M}^{0}$ and $M_2\in \mathcal{M}^{k}$, From this decomposition We deduce the functional equation
    $$M(x,y)=1+(2x+xy+x^2M(x,0))M(x,y),$$
    and so $$M(x,y)=\frac{2}{1-2x-2xy+\sqrt{1-4x}}.$$ 
   Finally, we obtain $I(x,y)$ and $I(x,1)$ with the shift $I(x,y)=1+xM(x,y)$.
\end{proof}



\subsection{In the lattice $\FF_n^2$}

According to Fact~\ref{fac:first-insertconv2}, the rules that describe the first descent lengths of paths $P'\le Q'$ obtained from  $[P,Q]\in\FF_{n}^\infty$ in terms of the first descent lengths $(a,b)$ of $P$ and $Q$, with respect to the difference $k$ of the lengths of ascents of $Q$ and $P$ is
$$\left\{\begin{array}{lll}
  (1,1)_k  \rightarrow  &(1,1)_0, (2,2)_0, (1,2)_1, &k=0,\\
  (1,1)_k  \rightarrow  &(1,1)_k, (2,2)_k, (1,2)_{k+1}, (2,1)_{k-1}, &k\geq 1,\\
  (1,2)_k  \rightarrow  &(1,1)_k, (2,1)_{k-1}, &k\geq 1,\\
   (2,1)_k  \rightarrow  &(1,1)_k, (1,2)_{k+1}, &k\geq 0,\\
    (2,2)_k  \rightarrow & (1,1)_k, &k\geq 0,
    \label{rule2}
    \end{array}\right.
  $$
starting with the root $(1,1)_0$ corresponding to the interval $[UD,UD]$. In fact, the above rules describe the construction of some  paths in the first quadrant of the plane.

\begin{theorem}\label{bij int motzkin fibo} There is a bijection between intervals in 
    $\FF^2_n$ and  bicolored Motzkin paths of length $n-1$ and avoiding the seven patterns $\Fg\Fg,\Fg\D,\Fg\U,\D\Fg,\U\Fg,\U\U,\D\D$.
\end{theorem}
\begin{proof}
Using the previous rules, we associate a bicolored Motzkin path with each generated interval by the following process. 
The bicolored Motzkin path is constructed from the origin $(0,0)$ using steps $\U=(1,1)$, $\D=(1,-1)$, $\Fb=(1,0)$, $\Fg=(1,0)$ corresponding, respectively,  to the rules $(a,b)\rightarrow (1,2)$, $(a,b)\rightarrow (2,1)$, $(a,b)\rightarrow (1,1)$ and $(a,b)\rightarrow (2,2)$. See Figure~\ref{ex bij int fibo motzkin} for an illustration of the bijection. Note that the subscript $k$ in $(a,b)_k$ corresponds to the final height of the Motzkin  path, and since it is always non-negative, the obtained lattice paths stay in the first quarter of the plane. Finally,  the one-to-one correspondence between interval and bicolored Motzkin paths is obtained with the avoidance of  the patterns $\Fg\Fg,\Fg\D,\Fg\U,\D\Fg,\U\Fg,\U\U,\D\D$.  Indeed, for example the above system does not have two consecutive rules of the form $(a,b)\rightarrow (2,2)$ which is equivalent to  the avoidance of $\Fg\Fg$. The other avoidances can be obtained {\it mutatis mutandis}.
\end{proof}

\begin{figure}[h]
    \centering
    \captionsetup{justification=centering,margin=0.3cm}
    \begin{tikzpicture}
    \node (1) at (0,0) {\begin{tikzpicture}[scale=0.3]
        \draw[solid,line width=0.4mm,blue] (0,0)-- ++\u -- ++\d;
        \draw[solid,line width=0.4mm,red] (0,0.6)-- ++\u -- ++\d;
        \filldraw[blue] (0,0)  circle (4pt);
        \filldraw[blue] (1,1)  circle (4pt);
        \filldraw[blue] (2,0)  circle (4pt);
        \filldraw[red] (0,0.5)  circle (4pt);
        \filldraw[red] (1,1.5)  circle (4pt);
        \filldraw[red] (2,0.5)  circle (4pt);
    \end{tikzpicture}};
    \node (2) at (2,0) {\begin{tikzpicture}[scale=0.3]
        \draw[solid,line width=0.4mm,blue] (0,0)-- ++\u -- ++\d-- ++\u -- ++\d;
        \draw[solid,line width=0.4mm,red] (0,0.6)-- ++\u-- ++\u -- ++\d -- ++\d;
        \filldraw[blue] (0,0)  circle (4pt);
        \filldraw[blue] (1,1)  circle (4pt);
        \filldraw[blue] (2,0)  circle (4pt);
        \filldraw[blue] (3,1)  circle (4pt);
        \filldraw[blue] (4,0)  circle (4pt);
        \filldraw[red] (0,0.5)  circle (4pt);
        \filldraw[red] (1,1.5)  circle (4pt);
        \filldraw[red] (2,2.5)  circle (4pt);
        \filldraw[red] (3,1.5)  circle (4pt);
        \filldraw[red] (4,0.5)  circle (4pt);
    \end{tikzpicture}};
    \node (3) at (4.6,0) {\begin{tikzpicture}[scale=0.3]
        \draw[solid,line width=0.4mm,blue] (0,0)-- ++\u -- ++\d-- ++\u -- ++\d-- ++\u -- ++\d;
        \draw[solid,line width=0.4mm,red] (0,0.6)-- ++\u-- ++\u -- ++\d-- ++\u  -- ++\d -- ++\d;
        \filldraw[blue] (0,0)  circle (4pt);
        \filldraw[blue] (1,1)  circle (4pt);
        \filldraw[blue] (2,0)  circle (4pt);
        \filldraw[blue] (3,1)  circle (4pt);
        \filldraw[blue] (4,0)  circle (4pt);
        \filldraw[blue] (5,1)  circle (4pt);
        \filldraw[blue] (6,0)  circle (4pt);
        \filldraw[red] (0,0.5)  circle (4pt);
        \filldraw[red] (1,1.5)  circle (4pt);
        \filldraw[red] (2,2.5)  circle (4pt);
        \filldraw[red] (3,1.5)  circle (4pt);
        \filldraw[red] (4,2.5)  circle (4pt);
        \filldraw[red] (5,1.5)  circle (4pt);
        \filldraw[red] (6,0.5)  circle (4pt);
    \end{tikzpicture}};
    \node (4) at (8,0) {\begin{tikzpicture}[scale=0.3]
        \draw[solid,line width=0.4mm,blue] (0,0)-- ++\u -- ++\u-- ++\d -- ++\d-- ++\u -- ++\d-- ++\u -- ++\d;
        \draw[solid,line width=0.4mm,red] (0,0.6)-- ++\u-- ++\u-- ++\u -- ++\d -- ++\d-- ++\u  -- ++\d -- ++\d;
        \filldraw[blue] (0,0)  circle (4pt);
        \filldraw[blue] (1,1)  circle (4pt);
        \filldraw[blue] (2,2)  circle (4pt);
        \filldraw[blue] (3,1)  circle (4pt);
        \filldraw[blue] (4,0)  circle (4pt);
        \filldraw[blue] (5,1)  circle (4pt);
        \filldraw[blue] (6,0)  circle (4pt);
        \filldraw[blue] (7,1)  circle (4pt);
        \filldraw[blue] (8,0)  circle (4pt);
        \filldraw[red] (0,0.5)  circle (4pt);
        \filldraw[red] (1,1.5)  circle (4pt);
        \filldraw[red] (2,2.5)  circle (4pt);
        \filldraw[red] (3,3.5)  circle (4pt);
        \filldraw[red] (4,2.5)  circle (4pt);
        \filldraw[red] (5,1.5)  circle (4pt);
        \filldraw[red] (6,2.5)  circle (4pt);
        \filldraw[red] (7,1.5)  circle (4pt);
        \filldraw[red] (8,0.5)  circle (4pt);
    \end{tikzpicture}};
    \node (5) at (12,0) {\begin{tikzpicture}[scale=0.3]
        \draw[solid,line width=0.4mm,blue] (0,0)-- ++\u-- ++\u -- ++\d -- ++\u-- ++\d -- ++\d-- ++\u -- ++\d-- ++\u -- ++\d;
        \draw[solid,line width=0.4mm,red] (0,0.6)-- ++\u-- ++\u-- ++\u -- ++\d-- ++\u -- ++\d -- ++\d-- ++\u  -- ++\d -- ++\d;
        \filldraw[blue] (0,0)  circle (4pt);
        \filldraw[blue] (1,1)  circle (4pt);
        \filldraw[blue] (2,2)  circle (4pt);
        \filldraw[blue] (3,1)  circle (4pt);
        \filldraw[blue] (4,2)  circle (4pt);
        \filldraw[blue] (5,1)  circle (4pt);
        \filldraw[blue] (6,0)  circle (4pt);
        \filldraw[blue] (7,1)  circle (4pt);
        \filldraw[blue] (8,0)  circle (4pt);
        \filldraw[blue] (9,1)  circle (4pt);
        \filldraw[blue] (10,0)  circle (4pt);
        \filldraw[red] (0,0.5)  circle (4pt);
        \filldraw[red] (1,1.5)  circle (4pt);
        \filldraw[red] (2,2.5)  circle (4pt);
        \filldraw[red] (3,3.5)  circle (4pt);
        \filldraw[red] (4,2.5)  circle (4pt);
        \filldraw[red] (5,3.5)  circle (4pt);
        \filldraw[red] (6,2.5)  circle (4pt);
        \filldraw[red] (7,1.5)  circle (4pt);
        \filldraw[red] (8,2.5)  circle (4pt);
        \filldraw[red] (9,1.5)  circle (4pt);
        \filldraw[red] (10,0.5)  circle (4pt);
    \end{tikzpicture}};
    \node (6) at (12.5,-2) {\begin{tikzpicture}[scale=0.275]
        \draw[solid,line width=0.4mm,blue] (0,0)-- ++\u-- ++\u-- ++\d -- ++\u -- ++\d -- ++\u-- ++\d -- ++\d-- ++\u -- ++\d-- ++\u -- ++\d;
        \draw[solid,line width=0.4mm,red] (0,0.6)-- ++\u-- ++\u-- ++\u-- ++\u -- ++\d -- ++\d-- ++\u -- ++\d -- ++\d-- ++\u  -- ++\d -- ++\d;
        \filldraw[blue] (0,0)  circle (4pt);
        \filldraw[blue] (1,1)  circle (4pt);
        \filldraw[blue] (2,2)  circle (4pt);
        \filldraw[blue] (3,1)  circle (4pt);
        \filldraw[blue] (4,2)  circle (4pt);
        \filldraw[blue] (5,1)  circle (4pt);
        \filldraw[blue] (6,2)  circle (4pt);
        \filldraw[blue] (7,1)  circle (4pt);
        \filldraw[blue] (8,0)  circle (4pt);
        \filldraw[blue] (9,1)  circle (4pt);
        \filldraw[blue] (10,0)  circle (4pt);
        \filldraw[blue] (11,1)  circle (4pt);
        \filldraw[blue] (12,0)  circle (4pt);
        \filldraw[red] (0,0.5)  circle (4pt);
        \filldraw[red] (1,1.5)  circle (4pt);
        \filldraw[red] (2,2.5)  circle (4pt);
        \filldraw[red] (3,3.5)  circle (4pt);
        \filldraw[red] (4,4.5)  circle (4pt);
        \filldraw[red] (5,3.5)  circle (4pt);
        \filldraw[red] (6,2.5)  circle (4pt);
        \filldraw[red] (7,3.5)  circle (4pt);
        \filldraw[red] (8,2.5)  circle (4pt);
        \filldraw[red] (9,1.5)  circle (4pt);
        \filldraw[red] (10,2.5)  circle (4pt);
        \filldraw[red] (11,1.5)  circle (4pt);
        \filldraw[red] (12,0.5)  circle (4pt);
    \end{tikzpicture}};
    \node (7) at (7,-2) {\begin{tikzpicture}[scale=0.275]
        \draw[solid,line width=0.4mm,blue] (0,0)-- ++\u-- ++\u-- ++\u -- ++\d-- ++\d -- ++\u -- ++\d -- ++\u-- ++\d -- ++\d-- ++\u -- ++\d-- ++\u -- ++\d;
        \draw[solid,line width=0.4mm,red] (0,0.6)-- ++\u-- ++\u-- ++\u-- ++\u-- ++\d -- ++\u -- ++\d -- ++\d-- ++\u -- ++\d -- ++\d-- ++\u  -- ++\d -- ++\d;
        \filldraw[blue] (0,0)  circle (4pt);
        \filldraw[blue] (1,1)  circle (4pt);
        \filldraw[blue] (2,2)  circle (4pt);
        \filldraw[blue] (3,3)  circle (4pt);
        \filldraw[blue] (4,2)  circle (4pt);
        \filldraw[blue] (5,1)  circle (4pt);
        \filldraw[blue] (6,2)  circle (4pt);
        \filldraw[blue] (7,1)  circle (4pt);
        \filldraw[blue] (8,2)  circle (4pt);
        \filldraw[blue] (9,1)  circle (4pt);
        \filldraw[blue] (10,0)  circle (4pt);
        \filldraw[blue] (11,1)  circle (4pt);
        \filldraw[blue] (12,0)  circle (4pt);
        \filldraw[blue] (13,1)  circle (4pt);
        \filldraw[blue] (14,0)  circle (4pt);
        \filldraw[red] (0,0.5)  circle (4pt);
        \filldraw[red] (1,1.5)  circle (4pt);
        \filldraw[red] (2,2.5)  circle (4pt);
        \filldraw[red] (3,3.5)  circle (4pt);
        \filldraw[red] (4,4.5)  circle (4pt);
        \filldraw[red] (5,3.5)  circle (4pt);
        \filldraw[red] (6,4.5)  circle (4pt);
        \filldraw[red] (7,3.5)  circle (4pt);
        \filldraw[red] (8,2.5)  circle (4pt);
        \filldraw[red] (9,3.5)  circle (4pt);
        \filldraw[red] (10,2.5)  circle (4pt);
        \filldraw[red] (11,1.5)  circle (4pt);
        \filldraw[red] (12,2.5)  circle (4pt);
        \filldraw[red] (13,1.5)  circle (4pt);
        \filldraw[red] (14,0.5)  circle (4pt);
    \end{tikzpicture}};
    \node (8) at (1,-2) {\begin{tikzpicture}[scale=0.275]
        \draw[solid,line width=0.4mm,blue] (0,0)-- ++\u-- ++\u-- ++\u -- ++\d-- ++\u -- ++\d-- ++\d -- ++\u -- ++\d -- ++\u-- ++\d -- ++\d-- ++\u -- ++\d-- ++\u -- ++\d;
        \draw[solid,line width=0.4mm,red] (0,0.6)-- ++\u-- ++\u-- ++\u-- ++\u-- ++\d -- ++\u-- ++\d -- ++\u -- ++\d -- ++\d-- ++\u -- ++\d -- ++\d-- ++\u  -- ++\d -- ++\d;
        \filldraw[blue] (0,0)  circle (4pt);
        \filldraw[blue] (1,1)  circle (4pt);
        \filldraw[blue] (2,2)  circle (4pt);
        \filldraw[blue] (3,3)  circle (4pt);
        \filldraw[blue] (4,2)  circle (4pt);
        \filldraw[blue] (5,3)  circle (4pt);
        \filldraw[blue] (6,2)  circle (4pt);
        \filldraw[blue] (7,1)  circle (4pt);
        \filldraw[blue] (8,2)  circle (4pt);
        \filldraw[blue] (9,1)  circle (4pt);
        \filldraw[blue] (10,2)  circle (4pt);
        \filldraw[blue] (11,1)  circle (4pt);
        \filldraw[blue] (12,0)  circle (4pt);
        \filldraw[blue] (13,1)  circle (4pt);
        \filldraw[blue] (14,0)  circle (4pt);
        \filldraw[blue] (15,1)  circle (4pt);
        \filldraw[blue] (16,0)  circle (4pt);
        \filldraw[red] (0,0.5)  circle (4pt);
        \filldraw[red] (1,1.5)  circle (4pt);
        \filldraw[red] (2,2.5)  circle (4pt);
        \filldraw[red] (3,3.5)  circle (4pt);
        \filldraw[red] (4,4.5)  circle (4pt);
         \filldraw[red] (5,3.5)  circle (4pt);
         \filldraw[red] (6,4.5)  circle (4pt);
         \filldraw[red] (7,3.5)  circle (4pt);
         \filldraw[red] (8,4.5)  circle (4pt);
         \filldraw[red] (9,3.5)  circle (4pt);
         \filldraw[red] (10,2.5)  circle (4pt);
         \filldraw[red] (11,3.5)  circle (4pt);
         \filldraw[red] (12,2.5)  circle (4pt);
         \filldraw[red] (13,1.5)  circle (4pt);
        \filldraw[red] (14,2.5)  circle (4pt);
         \filldraw[red] (15,1.5)  circle (4pt);
         \filldraw[red] (16,0.5)  circle (4pt);
    \end{tikzpicture}};
    \draw[->] (1) -- (2);
    \draw[->] (2) -- (3);
    \draw[->] (3) -- (4);
    \draw[->] (4) -- (5);
    \draw[->] (14,0) arc (90:-60:1);
    \draw[->] (6) -- (7);
    \draw[->] (7) -- (8);
    \node (a) at (0.9,0.5) {$\U$};
    \node (b) at (3.2,0.5) {$\Fg$};
    \node (c) at (6.2,0.5) {$\Fb$};
    \node (d) at (9.9,0.5) {$\Fg$};
    \node (e) at (14.5,-1) {$\U$};
    \node (f) at (9.95,-1.5) {$\D$};
    \node (g) at (4.2,-1.5) {$\Fg$};
    \end{tikzpicture}
    \begin{tikzpicture}
        \draw[solid,line width=0.4mm] (-0.5,0)--(7.5,0);
        \draw[solid,line width=0.4mm] (0,-0.5)--(0,2.5);
        \draw[dashed,line width=0.2mm] (0,1)--(7.3,1);
        \draw[dashed,line width=0.2mm] (0,2)--(7.3,2);
        \draw[dashed,line width=0.2mm] (1,0)--(1,2.3);
        \draw[dashed,line width=0.2mm] (2,0)--(2,2.3);
        \draw[dashed,line width=0.2mm] (3,0)--(3,2.3);
        \draw[dashed,line width=0.2mm] (4,0)--(4,2.3);
        \draw[dashed,line width=0.2mm] (5,0)--(5,2.3);
        \draw[dashed,line width=0.2mm] (6,0)--(6,2.3);
        \draw[dashed,line width=0.2mm] (7,0)--(7,2.3);
        \draw[solid,line width=0.4mm,blue,->] (0,0)--(1,1);
        \draw[solid,line width=0.4mm,Green,->] (1,1)--(2,1);
        \draw[solid,line width=0.4mm,blue,->] (2,1)--(3,1);
        \draw[solid,line width=0.4mm,Green,->] (3,1)--(4,1);
        \draw[solid,line width=0.4mm,blue,->] (4,1)--(5,2);
        \draw[solid,line width=0.4mm,blue,->] (5,2)--(6,1);
        \draw[solid,line width=0.4mm,Green,->] (6,1)--(7,1);
    \end{tikzpicture}
\caption{The generation of the interval $[U^2(UD)^2(DU)^2D^2(UD)^2,U^3(UD)^2(UD^2)^3]$ using the rules in the proof on Theorem~\ref{bij int motzkin fibo}. This interval is thus associated with the bicolored Motzkin path $\U\Fg\Fb\Fg\U\D\Fg$.}
\label{ex bij int fibo motzkin}
\end{figure}

\begin{thm}\label{int fibo}
    The generating function $J(x,y)$ for the number of intervals $[P,Q]$ in $\F^2_n$ with respect to $n$ and the difference of the lengths of ascents of $Q$ and $P$ is given by
$$J(x,y)=\frac{1+x-x^2+\sqrt{x^4-2x^3-x^2-2x+1}}{(1 - x^2)\sqrt{x^4 - 2x^3 - x^2 - 2x + 1} + 1 + x^4 - x^3 - 2(y + 1)x^2 - x}.$$

    The generating function $J(x,1)$ for the number of intervals $[P,Q]$ in $\F^2_n$ is
    $$J(x,1)=\frac{-x^2+3x-1+\sqrt{x^4 - 2x^3 - x^2 - 2x + 1}}{2x(x^2 - 3x + 1)(x + 1)}.$$
    The coefficient of $x^n$ in the series expansion is asymptotically 
    
    $$\frac{11+5\sqrt{5}}{20}\sqrt{\frac{14\sqrt{5}-30}{\pi}}\cdot n^{-1/2}\left(\frac{3+\sqrt{5}}{2}\right)^n.$$

\end{thm}

\begin{proof}
    By Theorem~\ref{bij int motzkin fibo}, it suffices to count bicolored Motzkin paths avoiding the patterns $\Fg\Fg,\Fg\D,\Fg\U,\D\Fg,\U\Fg,\U\U,\D\D$ with respect to the number of steps and the final height. 
    We classify those paths in 9 different categories:
    \begin{enumerate}[label=(\roman*)]
        \item the empty path,
        \item paths starting with an $\Fb$, followed by any path of the class,
        \item the path $\Fg$,
        \item paths starting with $\Fg\Fb$, followed by any path of the class,
        \item the path $\U$,
        \item paths starting with $\U\D$, followed by a path of the class not starting with $\Fg$,
        \item paths starting with $\U\Fb$, followed by any path of the class,
        \item paths starting with $\U\Fb\D$, followed by a path of the class not starting with $\Fg$,
        \item paths starting with $\U\Fb$, followed by any path of the class ending at height 0, followed by $\Fb\D$, followed by any path of the class not starting with $\Fg$.
    \end{enumerate}
    This decomposition gives the following equation:
    \begin{align*}
        A(x,y)=1&+xA(x,y)+x+x^2A(x,y)+xy+x^2(A(x,y)-x-x^2A(x,y))\\
        &+x^2yA(x,y)+x^3(A(x,y)-x-x^2A(x,y))\\
        &+x^4A(x,0)(A(x,y)-x-x^2A(x,y)).
    \end{align*}
    By specializing $y=0$ and solving this equation, we find 
    $$A(x,0)=\frac{1-x-x^2-2x^3-\sqrt{x^4-2x^3-x^2-2x+1}}{2x^4}.$$
    By plugging this expression in the previous equation we can then solve for $A(x,y)$, which gives
    $$A(x,y)=\frac{x\sqrt{x^4-2x^3-x^2-2x+1} + 2 - x^3 + x^2 + (2y + 1)x}{(1 - x^2)\sqrt{x^4 - 2x^3 - x^2 - 2x + 1} + 1 + x^4 - x^3 - 2(y + 1)x^2 - x}.$$
    Then we obtain $J(x,y)$ from $A(x,y)$ with the shift $J(x,y)=1+xA(x,y)$.
\end{proof}

The first terms of the series expansion of $J(x,1)$ are $$1+x+3x^2+6x^3+15x^4+35x^5+86x^6+210x^7+520x^8+1292x^9+O(x^{10}),$$ and the sequence of coefficients does not appear in \cite{oeis}.



\subsection{In the lattice $\FF_n^p, p\geq 3$}

The following theorem gives a generalization of  Theorem~\ref{bij int motzkin fibo}, providing a bijection between intervals in $\FF_n^p$ for any $p\geq 3$, and some bicolored Motzkin paths. However, obtaining a count of those paths seems challenging, because of the exponential number of forbidden patterns.

\begin{theorem}\label{bij int motzkin fibop3} There is a bijection between intervals in 
    $\FF^p_n$ and  bicolored Motzkin paths of length $n-1$ and avoiding the $2^{p+1}-1$ patterns of the set $\{\Fg,\U\}^p\cup\{\Fg,\D\}^p$.
\end{theorem}
\begin{proof}
    An interval $[P,Q]$ in $\FF_n^p$ can be generated with similar rules as in Theorem \ref{cor:bij-mn}. However we need to be careful not creating an occurrence of $D^{p+1}$ in $P$ or $Q$. Such an occurrence appears in $P$ if and only if, with the notation of Fact \ref{fac:first-insertconv1}, we insert $p$ times consecutively the factor $UD$ at height $i$ (and similarly with $Q$ at height $j$). The steps corresponding to an insertion at height $i$ in $P$ (resp. $j$ in $Q$) are $\Fg$ and $\D$ (resp. $\Fg$ and $\U$). The theorem then follows.
\end{proof}

\section{Bijections with other combinatorial objects}\label{sec:bijections}

The generalized Fibonacci numbers count a variety of combinatorial objects, so with no surprise, there are a lot of bijections between $\F_n^p$ and other objects. In this section we present some of them, with well-known objects, and we show how the lattice structure is conveyed. As a consequence, Theorems \ref{int binary} and \ref{int fibo} also give the enumeration of the intervals viewed as in  Propositions \ref{int catalan}, \ref{int compositions} and \ref{int powerset}, for $p\in\{2,\infty\}$.

\subsection{Catalan words} 
\begin{figure}[h]
    \centering
    \captionsetup{justification=centering,margin=1cm}
    \begin{tikzpicture}[scale=0.6]
        \draw[solid,line width=0.4mm] (0,0)-- ++\u -- ++\u -- ++\u -- ++\u -- ++\u -- ++\d -- ++\u -- ++\d -- ++\d -- ++\d -- ++\u -- ++\d -- ++\d -- ++\d;
        \filldraw (0,0)  circle (4pt);
        \filldraw (1,1)  circle (4pt);
        \filldraw (2,2)  circle (4pt);
        \filldraw (3,3)  circle (4pt);
        \filldraw (4,4)  circle (4pt);
        \filldraw (5,5)  circle (4pt);
        \filldraw (6,4)  circle (4pt);
        \filldraw (7,5)  circle (4pt);
        \filldraw (8,4)  circle (4pt);
        \filldraw (9,3)  circle (4pt);
        \filldraw (10,2)  circle (4pt);
        \filldraw (11,3)  circle (4pt);
        \filldraw (12,2)  circle (4pt);
        \filldraw (13,1)  circle (4pt);
        \filldraw (14,0)  circle (4pt);
        \node[blue] at (5.8,4.8) {$2$};
        \node[blue] at (7.8,4.8) {$1$};
        \node[blue] at (8.8,3.8) {$1$};
        \node[blue] at (9.8,2.8) {$1$};
        \node[blue] at (11.8,2.8) {$0$};
        \node[blue] at (12.8,1.8) {$0$};
        \node[blue] at (13.8,0.8) {$0$};
        \draw[dashed] (0,0)--(14,0);
    \end{tikzpicture}
    \caption{The path $P=U^5D(UD^3)^2\in\F_7^\infty$ is associated with the Catalan word $w(P)=0001112$.}
    \label{bijection paths words}
\end{figure}
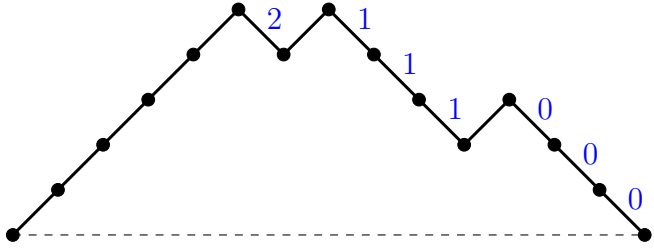

A length $n$ \textit{Catalan word} is a word $w_1\ldots w_n$ over the set of non-negative integers, with $w_1=0$ and $0\leq w_i\leq w_{i-1}+1$ for $i=2,3,\ldots,n$. We present a bijection between $\F_n^p$ and the set $\mathcal{C}_n^p$ of length $n$ non-decreasing Catalan words avoiding $p+1$ consecutive occurrences of the same letter. Similarly, $\F_n^\infty$ is in bijection with the set $\mathcal{C}_n^\infty$ of length $n$ and non-decreasing Catalan words. Let $P$ be a Dyck path in $\F_n^p$ . Label the $n$ down steps of $P$ from $1$ to $n$ and from right to left. For $i=1,\ldots,n$, let $w_i$ be the number of up steps in $P$ that are at the right of the down step number $i$. We set $w(P)=w_1\ldots w_n$. Since $P$ avoids $DUU$, $w(P)$ is a Catalan word. It is clearly non-decreasing by construction, and since $P$ avoids $D^{p+1}$, $w(P)$ does not have $p+1$ consecutive occurrences of the same letter. Conversely, let $w=w_1\ldots w_n\in\mathcal{C}_n^p$. Let $k:=w_n$ be the greatest letter in $w$, and for $i=0,\ldots,k$, let $a_i$ be the number of occurrences of the letter $i$ in $w$. Let $P=U^{n-k}D^{a_k}UD^{a_{k-1}}\ldots UD^{a_1}UD^{a_0}$. It is easy to check that $w(P)=w$, $P$ avoids $DUU$, and since $a_i\leq p$ for $0\leq i\leq k$, $P$ avoids $D^{p+1}$, thus $P\in\F_n^p$. Then, $w$ is a bijection between $\F_n^p$ and $\C_n^p$ (for $p\geq2$ and $p=\infty$), see Figure \ref{bijection paths words}.

As a consequence, $\FF_n^p$ induces a lattice structure on $\mathcal{C}_n^p$, the cover relation being 
$$v\lessdot w \Longleftrightarrow \mbox{ there exists } i \mbox{ such that } v_i=w_i+1, \mbox{ and } v_j=w_j \mbox{ for }j\neq i.$$
See Figure \ref{lattices on power sets} for an illustration. We then deduce the following proposition.

\begin{proposition}\label{int catalan}
    Let $v,w\in\mathcal{C}_n^p$. Then $[v,w]$ is an interval in $\FF_n^p$ if and only if for all $1\leq i \leq n$, $w_i\leq v_i$.
\end{proposition}

\subsection{Compositions} 
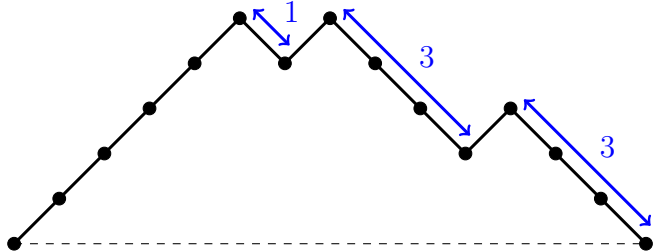
\begin{figure}[h]
    \centering
    \captionsetup{justification=centering,margin=1cm}
    \begin{tikzpicture}[scale=0.6]
        \draw[solid,line width=0.4mm] (0,0)-- ++\u -- ++\u -- ++\u -- ++\u -- ++\u -- ++\d -- ++\u -- ++\d -- ++\d -- ++\d -- ++\u -- ++\d -- ++\d -- ++\d;
        \filldraw (0,0)  circle (4pt);
        \filldraw (1,1)  circle (4pt);
        \filldraw (2,2)  circle (4pt);
        \filldraw (3,3)  circle (4pt);
        \filldraw (4,4)  circle (4pt);
        \filldraw (5,5)  circle (4pt);
        \filldraw (6,4)  circle (4pt);
        \filldraw (7,5)  circle (4pt);
        \filldraw (8,4)  circle (4pt);
        \filldraw (9,3)  circle (4pt);
        \filldraw (10,2)  circle (4pt);
        \filldraw (11,3)  circle (4pt);
        \filldraw (12,2)  circle (4pt);
        \filldraw (13,1)  circle (4pt);
        \filldraw (14,0)  circle (4pt);
        \draw[blue,<->,line width=0.4mm] (5.3,5.2)--(6.1,4.4);
        \draw[blue,<->,line width=0.4mm] (7.3,5.2)--(10.1,2.4);
        \draw[blue,<->,line width=0.4mm] (11.3,3.2)--(14.1,0.4);
        \node[blue] at (6.15,5.15) {$1$};
        \node[blue] at (9.15,4.15) {$3$};
        \node[blue] at (13.15,2.15) {$3$};
        \draw[dashed] (0,0)--(14,0);
    \end{tikzpicture}
    \caption{The path $P=U^5D(UD^3)^2\in\F_7^\infty$ is associated with the composition $\lambda(P)=(3,3,1)$.}
    \label{bijection paths composition}
\end{figure}

In this section we present a bijection between the elements of $\F_n^p$ and the compositions of $n$ with parts in $[1,p]$, and between $\F_n^\infty$ and all compositions of $n$. This one is quite natural, since any path $P\in\F_n^p$ can be uniquely written $P=U^{n-k+1}D^{\lambda_k}UD^{\lambda_{k-1}}\ldots UD^{\lambda_1}$ with $\lambda_1,\ldots,\lambda_k\in [1,p]$, and so $\lambda(P):=(\lambda_1,\ldots,\lambda_k)$ is a composition of $n$ with parts in $[1,p]$. We then easily see that $P\mapsto \lambda(P)$ is a bijection, see Figure \ref{bijection paths composition}. Thus $\FF_n^p$ induces a lattice structure on the compositions of $n$ with parts in $[1,p]$, the cover relation being 
$$(\lambda_1,\ldots,\lambda_k)\lessdot(\mu_1,\ldots,\mu_\ell)\Longleftrightarrow 
\left\{\begin{array}{l}
k=\ell,  \text{ and there exists } i\in [2,k] \ \text{ such that } \\ \quad \qquad \lambda_i>1, \ (\mu_{i-1},\mu_i)=(\lambda_{i-1}+1,\lambda_{i}-1), \\ \quad \qquad \mbox{and } \ \mu_j=\lambda_j \ \mbox{ for } j\not\in\{i-1,i\}, \mbox{ or }\\
 \ell=k-1,\ \lambda_k=1, \ \mu_\ell=\lambda_{k-1}+1, \\
 \qquad\quad \mbox{and } \mu_j=\lambda_j \ \mbox{ for } \ j\in [1,k-2].
 \end{array}\right.$$
See Figure \ref{lattices on power sets} for an illustration. From this we deduce the following proposition.

\begin{proposition}\label{int compositions}
    The order induced by $\FF_n^p$ on the compositions of $n$ with parts in $[1,p]$ is known as the \textit{dominance} order \cite{blass sagan}, defined by $\lambda\leq\mu$ if and only if for all $k$ we have $\sum_{i=1}^k\lambda_i\leq\sum_{i=1}^k\mu_i$.
\end{proposition}

\subsection{Power set of $[1,n-1]$} 
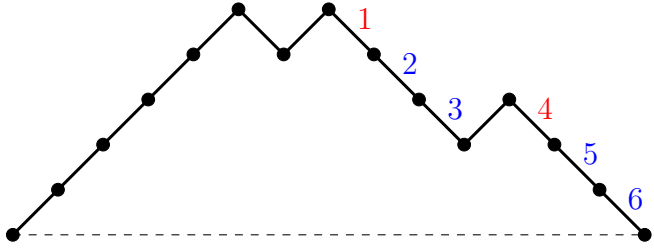
\begin{figure}[h]
    \centering
    \captionsetup{justification=centering,margin=1cm}
    \begin{tikzpicture}[scale=0.6]
        \draw[solid,line width=0.4mm] (0,0)-- ++\u -- ++\u -- ++\u -- ++\u -- ++\u -- ++\d -- ++\u -- ++\d -- ++\d -- ++\d -- ++\u -- ++\d -- ++\d -- ++\d;
        \filldraw (0,0)  circle (4pt);
        \filldraw (1,1)  circle (4pt);
        \filldraw (2,2)  circle (4pt);
        \filldraw (3,3)  circle (4pt);
        \filldraw (4,4)  circle (4pt);
        \filldraw (5,5)  circle (4pt);
        \filldraw (6,4)  circle (4pt);
        \filldraw (7,5)  circle (4pt);
        \filldraw (8,4)  circle (4pt);
        \filldraw (9,3)  circle (4pt);
        \filldraw (10,2)  circle (4pt);
        \filldraw (11,3)  circle (4pt);
        \filldraw (12,2)  circle (4pt);
        \filldraw (13,1)  circle (4pt);
        \filldraw (14,0)  circle (4pt);
        \node[red] at (7.8,4.8) {$1$};
        \node[blue] at (8.8,3.8) {$2$};
        \node[blue] at (9.8,2.8) {$3$};
        \node[red] at (11.8,2.8) {$4$};
        \node[blue] at (12.8,1.8) {$5$};
        \node[blue] at (13.8,0.8) {$6$};
        \draw[dashed] (0,0)--(14,0);
    \end{tikzpicture}
    \caption{The path $P=U^5D(UD^3)^2\in\F_7^\infty$ is associated with the subset $A(P)=\{2,3,5,6\}\subseteq\{1,\ldots,6\}$.}
    \label{bijection paths subsets}
\end{figure}

Here we present a bijection between $\F_n^p$, $p\in\N_{\geq2}\cup\{\infty\}$, and the subsets of $[1,n-1]$ having no $p$ consecutive elements (when $p=\infty$ it is just the whole powerset). Let $P\in\F_n^p$ be a Dyck path. We write $P=U^iDQ$ with $i$ the length of the first run of $U$'s in $P$, and $Q$ a path with $n-1$ $D$'s and $n-i$ $U$'s, avoiding $UU$. Then, label the $D$'s of $Q$ from 1 to $n-1$ and from left to right. Let $A(P)\subseteq[1,n-1]$ be the set of labels of the $D$'s that are not preceded by a $U$, see Figure \ref{bijection paths subsets} for an example. It is easy to check that if $P\in\F_n^p$ then $A(P)$ does not have $p$ consecutive elements.

\begin{figure}[h]
\centering
\captionsetup{justification=centering,margin=1cm}
\begin{tikzpicture}[scale=1]
    \node (1234) at (0,10) {$\{1,2,3,4\}$};
    \node (234) at (0,9) {$\{2,3,4\}$}; 
    \node (134) at (0,8) {$\{1,3,4\}$};
    \node (124) at (-0.866,7) {$\{1,2,4\}$};
    \node (34) at (0.866,7) {$\{3,4\}$};
    \node (123) at (-0.866,6) {$\{1,2,3\}$};
    \node (24) at (0.866,6) {$\{2,4\}$};
    \node (23) at (-0.866,5) {$\{2,3\}$};
    \node (14) at (0.866,5) {$\{1,4\}$};
    \node (13) at (-0.866,4) {$\{1,3\}$};
    \node (4) at (0.866,4) {$\{4\}$};
    \node (12) at (-0.866,3) {$\{1,2\}$};
    \node (3) at (0.866,3) {$\{3\}$};
    \node (2) at (0,2) {$\{2\}$};
    \node (1) at (0,1) {$\{1\}$};
    \node (0) at (0,0) {$\emptyset$};
    \draw (1234) -- (234);\draw (234) -- (134);\draw (134) -- (124);\draw (134) -- (34);\draw (124) -- (123);\draw (34) -- (24);\draw (124) -- (24);\draw (123) -- (23);\draw (24) -- (14);\draw (24) -- (23);\draw (23) -- (13);\draw (14) -- (4);\draw (14) -- (13);\draw (13) -- (12);\draw (4) -- (3);\draw (13) -- (3);\draw (3) -- (2);\draw (12) -- (2);\draw (2) -- (1);\draw (1) -- (0);
\end{tikzpicture}
\hspace{2cm}
\begin{tikzpicture}[scale=1]
    \node (1234) at (0,10) {$00000$};
    \node (234) at (0,9) {$00001$}; 
    \node (134) at (0,8) {$00011$};
    \node (124) at (-0.866,7) {$00111$};
    \node (34) at (0.866,7) {$00012$};
    \node (123) at (-0.866,6) {$01111$};
    \node (24) at (0.866,6) {$00112$};
    \node (23) at (-0.866,5) {$01112$};
    \node (14) at (0.866,5) {$00122$};
    \node (13) at (-0.866,4) {$01122$};
    \node (4) at (0.866,4) {$00123$};
    \node (12) at (-0.866,3) {$01222$};
    \node (3) at (0.866,3) {$01123$};
    \node (2) at (0,2) {$01223$};
    \node (1) at (0,1) {$01233$};
    \node (0) at (0,0) {$01234$};
    \draw (1234) -- (234);\draw (234) -- (134);\draw (134) -- (124);\draw (134) -- (34);\draw (124) -- (123);\draw (34) -- (24);\draw (124) -- (24);\draw (123) -- (23);\draw (24) -- (14);\draw (24) -- (23);\draw (23) -- (13);\draw (14) -- (4);\draw (14) -- (13);\draw (13) -- (12);\draw (4) -- (3);\draw (13) -- (3);\draw (3) -- (2);\draw (12) -- (2);\draw (2) -- (1);\draw (1) -- (0);
\end{tikzpicture}
\hspace{2cm}
\begin{tikzpicture}[scale=1]
    \node (1234) at (0,10) {$(5)$};
    \node (234) at (0,9) {$(4,1)$}; 
    \node (134) at (0,8) {$(3,2)$};
    \node (124) at (-0.866,7) {$(2,3)$};
    \node (34) at (0.866,7) {$(3,1,1)$};
    \node (123) at (-0.866,6) {$(1,4)$};
    \node (24) at (0.866,6) {$(2,2,1)$};
    \node (23) at (-0.866,5) {$(1,3,1)$};
    \node (14) at (0.866,5) {$(2,1,2)$};
    \node (13) at (-0.866,4) {$(1,2,2)$};
    \node (4) at (0.866,4) {$(2,1,1,1)$};
    \node (12) at (-0.866,3) {$(1,1,3)$};
    \node (3) at (0.866,3) {$(1,2,1,1)$};
    \node (2) at (0,2) {$(1,1,2,1)$};
    \node (1) at (0,1) {$(1,1,1,2)$};
    \node (0) at (0,0) {$(1,1,1,1,1)$};
    \draw (1234) -- (234);\draw (234) -- (134);\draw (134) -- (124);\draw (134) -- (34);\draw (124) -- (123);\draw (34) -- (24);\draw (124) -- (24);\draw (123) -- (23);\draw (24) -- (14);\draw (24) -- (23);\draw (23) -- (13);\draw (14) -- (4);\draw (14) -- (13);\draw (13) -- (12);\draw (4) -- (3);\draw (13) -- (3);\draw (3) -- (2);\draw (12) -- (2);\draw (2) -- (1);\draw (1) -- (0);
\end{tikzpicture}
\caption{The lattice $\FF_5^\infty$ on the power set of $\{1,2,3,4\}$ (left), the non-decreasing Catalan words of length $5$ (center) and the compositions of $5$ (right).}
\label{lattices on power sets}
\end{figure}
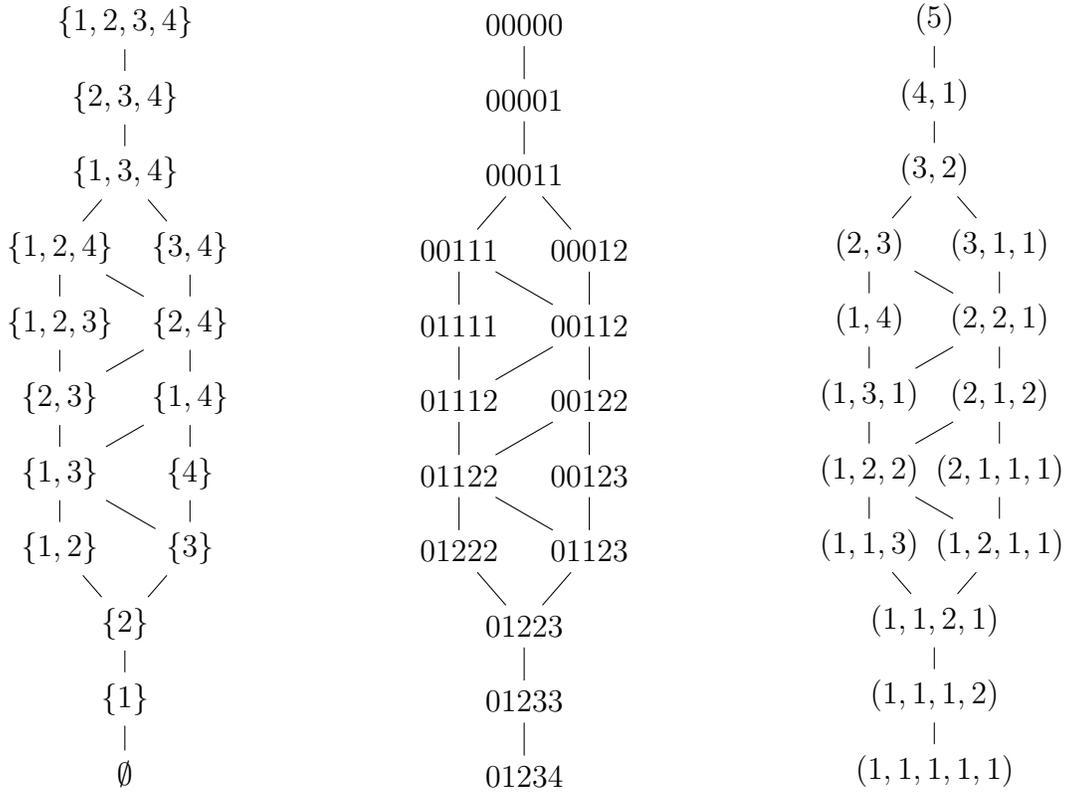

Conversely, if $A\subseteq[1,n-1]$, we build the following path: we start with $U^{|A|+1}D$, and then we add a path $Q=Q_1Q_2\ldots Q_{n-1}$ so that $Q_i=D$ if $i\in A$, and $Q_i=UD$ if $i\not\in A$. If $A$ does not contain $p$ consecutive elements, then it is easy to check that this path belongs to $\F_n^p$. We then obtain a lattice on the subsets of $[1,n-1]$ not having $p$ consecutive elements (all the subsets when $p=\infty$) with the following covering relation (see Figure \ref{lattices on power sets}): 

$$A\lessdot B \Longleftrightarrow
\left\{\begin{array}{l}
1\not\in A \mbox{ and } B=A\cup\{1\}, \mbox{ or }\\
\mbox{there exists a unique } x\in A \mbox{ such that } B=\{x+1\}\cup A\backslash\{x\}.\end{array}\right.$$

\begin{remark}\label{rank}
    In this setting, it is easy to see that the lattice is graded, the rank function $\rho$ being defined by $\rho(A)=\sum_{x\in A}x$. When $p=\infty$, the maximal element is $\{1,2,\ldots,n-1\}$, so the rank of $\FF_n^\infty$ is $n(n-1)/2$. When $p\in\N_{\geq2}$, the maximal element is $\{1,2,\ldots,n-1\}\backslash\{n-p,n-2p,\ldots,n-\left\lfloor\frac{n-1}{p}\right\rfloor p\}$, so the rank of $\FF_n^p$ is $$\rho(\FF_n^p)=\frac{n(n-1)}{2}-\left\lfloor\frac{n-1}{p}\right\rfloor\cdot\left(n-\frac{p\left(\left\lfloor\frac{n-1}{p}\right\rfloor+1\right)}{2}\right).$$
\end{remark}

For a set $A\subseteq [1,n-1]$, we denote by $A^c=[1,n-1]\backslash A$ its complement.

\begin{theorem}\label{symmetry F_n^infty}
    For all $A,B\in\FF_n^\infty$, $A\lessdot B$ if and only if $B^c\lessdot A^c$. So complement is a reverse ordering involution on $\FF_n^\infty$.
\end{theorem}
\begin{proof}
    $B=A\cup\{1\}$ if and only if $A^c=B^c\cup\{1\}$, and $B=\{x+1\}\cup A\backslash\{x\}$ if and only if $A^c=\{x+1\}\cup B^c\backslash\{x\}$.
\end{proof}

\begin{proposition}\label{int powerset}
    If $A,B\subseteq[1,n-1]$, then $[A,B]$ is an interval in $\FF_n^p$ if and only if $A=B$ or
    \begin{enumerate}[label=(\roman*)]
        \item $\sum_{x\in A}x<\sum_{x\in B}x$,
        \item $|A|\leq |B|$,
        \item $|A\backslash\{1\}|\leq |B\backslash\{1\}|$, and
        \item if $A=\{x_1>x_2>\cdots>x_{|A|}\}$ and $B=\{y_1>y_2>\cdots>y_{|B|}\}$, then $(x_1,\ldots,x_{|A|})\leq_{\text{lex}}(y_1,\ldots,y_{|A|})$, where $\leq_{\text{lex}}$ is the lexicographic order.
    \end{enumerate}
\end{proposition}
\begin{proof}
    Let $P,Q\in\F_n^p$ such that $[P,Q]$ is an interval in $\FF_n^p$, with $P\ne Q$. Let us prove that $[A(P),A(Q)]$ satisfies \textit{(i-iv)}. Since $A(P)\ne A(Q)$, we necessarily have $\rho(A(P))<\rho(A(Q))$, hence \textit{(i)}. Since $P$ lies below $Q$, the length of its first run of $U$'s is lower that the one of $Q$, hence \textit{(ii)}. Now that we have \textit{(ii)}, the only possibility for \textit{(iii)} to be false is that $|A(P)|=|A(Q)|$, $1\in A(Q)$ and $1\not\in A(P)$. But this would mean that $P$ (resp. $Q$) has the prefix $U^kDU$ (resp. $U^kD^2$) for some $k$, a contradiction to $P\leq Q$. Let $A(P)=\{x_1>\ldots>x_k\}$ and $A(Q)=\{y_1>\ldots>y_l\}$, with $l\geq k$. Let $i$ be such that $x_1=y_1,\ldots,x_{i-1}=y_{i-1}$ and $x_i\ne y_i$. Suppose that $x_i>y_i$. This means that there exists a path $S$ such that $P$ has the suffix $D^2S$, and $Q$ has the suffix $UDS$, again a contradiction to $P\leq Q$, hence \textit{(iv)}. Conversely, with similar arguments as before, we indeed have that if $A(P)$ and $A(Q)$ satisfy \textit{(i-iv)}, then $P<Q$.
\end{proof}

\begin{thebibliography}{99}

\bibitem{Aig} M. Aigner. 
\newblock Tur{\'a}n's graph theorem.
\newblock {\em Am. Math. Monthly}, 102 (1995), 808--816.

\bibitem{Barc} E. Barcucci, A. Bernini, L. Ferrari, and M. Poneti. 
\newblock A distributive lattice structure connecting Dyck paths, noncrossing partitions and $312$-avoiding permutations. \newblock {\em Order}, 22(4) (2005), 311--328.

\bibitem{Barp1} J.-L. Baril, and J.-M. Pallo.
\newblock The Phagocyte Lattice of Dyck Words.
\newblock {\em Order}, 23(2-3) (2006), 97--107.

\bibitem{Barp2} J.-L. Baril, and J.-M. Pallo.
\newblock The pruning-grafting lattice of binary trees.
\newblock {\em Theoretical Computer Science},  409(3) (2008),  382--293.

\bibitem{BKN}  J.-L. Baril, S. Kirgizov, and M. Naima. 
\newblock A lattice on Dyck paths close to the Tamari lattice.
\newblock \href{https://arxiv.org/abs/2309.00426}{https://arxiv.org/abs/2309.00426}, 2024.

\bibitem{BBKN} J.-L. Baril, M. Bousquet-M{\'e}lou, S. Kirgizov, and M. Naïma.
\newblock The ascent lattice on Dyck paths.
\newblock \href{https://arxiv.org/abs/2409.15982}{https://arxiv.org/abs/2409.15982}, 2024.

\bibitem{Barmo} J.-L. Baril, and J.-M. Pallo.
\newblock A Motzkin filter in the Tamari lattice.
\newblock {\em Discrete mathematics}, 338 (2015), 1370--1378.

\bibitem{Berg} F. Bergeron, and L.-F. Pr{\'e}ville-Ratelle.
\newblock Higher trivariate diagonal harmonics via generalized Tamari posets.
\newblock {\em J. comb.} 3(3) (2012), 317--341.

\bibitem{Bern} O. Bernardi, and N. Bonichon.
\newblock Intervals in Catalan lattices and realizers of triangulations.
\newblock {\em Journal of Combinatorial Theory, Series A}, 116 (2009), 55–-75.

\bibitem{Vaj} A. Bernini, S. Bilotta, R. Pinzani, V. Vajnovszki.
\newblock A trace partitioned Gray code for $q$-ary generalized Fibonacci strings.
\newblock {\em Discrete Mathematical Sciences and Cryptography}, 18(6) (2015) 751--761. 

\bibitem{blass sagan} A. Blass, B.E. Sagan.
\newblock M{\"o}bius Functions of Lattices.
\newblock {\em Advances in Mathematics}, 127, (1997), 94--123.

\bibitem{bollobas} B. Bollob\'as, \emph{Extremal graph theory}, Academic Press, 1978.

\bibitem{Bous} M. Bousquet-M{\'e}lou, E. Fusy, and L.-F. Pr{\'e}ville-Ratelle.
\newblock The number of intervals in the $m$-Tamari lattices.
\newblock {\em Electron. J. Combin.}, 18(2) (2011), paper 31.

\bibitem{Boucha} M. Bousquet-M{\'e}lou, and F. Chapoton.
\newblock Intervals in the greedy Tamari posets.
\newblock \href{https://arxiv.org/abs/2303.18077}{https://arxiv.org/abs/2303.18077}, 2024.

\bibitem{BFT} M. Bouvel, L. Ferrari, and B.E. Tenner.
\newblock Between weak and Bruhat: middle order on permutations.
\newblock \href{https://arxiv.org/abs/2405.08943}{https://arxiv.org/abs/2405.08943}, 2024.


\bibitem{Chap} F. Chapoton.
\newblock Sur le nombre d’intervalles dans les treillis de Tamari.
\newblock {\em S{\'e}m. Lothar. Combin.}, 55 (2006), Art. B55f (electronic).

\bibitem{Chapo} F. Chapoton.
\newblock Some properties of a new partial order on Dyck paths.
\newblock {\em Algebraic combinatorics}, 3(2) (2020), 433--463.

\bibitem{chen1} C. Chenevi{\`e}re.
\newblock Linear intervals in the Tamari and the Dyck lattices and in the alt-Tamari posets.
\newblock \href{https://arxiv.org/abs/2209.00418}{https://arxiv.org/abs/2209.00418}, 2022.

\bibitem{chen2} C. Chenevi{\`e}re. Enumerative study of intervals in lattices of Tamari type. PhD thesis, IRMA, Universit{\'e} de
Strasbourg, 2023. 

\bibitem{Fan} W. Fang, and L.-F. Pr{\'e}ville-Ratelle.
\newblock The enumeration of generalized Tamari intervals.
\newblock {\em European J. Combin.}, 61 (2017), 69--84.

\bibitem{analytic combinatorics}
P. Flajolet and R. Sedgewick, \emph{Analytic Combinatorics}. Cambridge University Press, 2005.

\bibitem{Grat} G. Gr{\"a}tzer. {\em General Lattice Theory}, Second edition, Birkh{\"a}user, 1998.

\bibitem{Huan} S. Huang, and D. Tamari.
\newblock Problems of associativity: A simple proof for the lattice property of systems ordered by a semi-associative law.
\newblock {\em J. Combinatorial Theory Ser. A}, 13 (1972), 7--13.

\bibitem{Knu} D.E. Knuth, The Art of Computer Programming, Fundamental Algorithms, Vol. I. Reading, Mass.: Addison-Wesley, 1973.

\bibitem{Kos} T. Koshy, Fibonacci and Lucas Numbers with Applications, A Wiley-Interscience Publication, 2001.

\bibitem{KO} D. Kremer, K. O’Hara.
\newblock A bijection between maximal chains in Fibonacci posets.
\newblock {\em J. Combin. Theory(Series A)}, 78(2) (1997) 268-–279.


\bibitem{Kre} D. Kremer.
\newblock A bijection between intervals in the Fibonacci posets.
\newblock {\em Discrete Math.},  217 (2000) 225–-235.

\bibitem{Mil} E. Miles.
\newblock Generalized Fibonacci numbers and associated matrices.
\newblock {\em The American
Mathematical Monthly}, 67(8), (1960) 745–-752.

\bibitem{Simi} R. Simion, and D. Ullman.
\newblock On the structure of lattice of noncrossing partitions.
\newblock {\em Discrete Math.}, 98 (1991), 193--206

\bibitem{oeis} N.J.A. Sloane, OEIS Foundation Inc., The On-line Encyclopedia of Integer Sequences, available electronically at \href{http://oeis.org}{http://oeis.org}.

\bibitem{Stan} R.P. Stanley.
\newblock The Fibonacci Lattice.
\newblock {\em Fibonacci Quart.}, 13 (1975) 215–-232.

\bibitem{Stanley} R.P. Stanley.
\newblock Differential posets.
\newblock {\em Journal of the American Mathematical Society}, 1(4) (1988), 919–-961, 


\bibitem{enum comb} R.P. Stanley, \emph{Enumerative Combinatorics, Volume 1}, Cambridge Studies in Advanced Mathematics 49, Cambridge University Press, Cambridge, 2012.

\bibitem{Stfur} R.P. Stanley.
\newblock Further Combinatorial Properties of Two Fibonacci Lattices.
\newblock {\em European Journal of Combinatorics}, 11(2) (1990), 181--188.

\bibitem{Tam} D. Tamari.
\newblock The algebra of bracketings and their enumeration.
\newblock {\em Nieuw Archief voor Wiskunde}, 10 (1962), 131--146.

\bibitem{turan} P. Tur\'an.
\newblock On an extremal problem in graph theory (in Hungarian).
\newblock {\em  Math. Fiz. Lapok}, 48 (1941), 436–-452.





\end{thebibliography}
\end{document}